\documentclass[twoside,a4paper,12pt]{article}
\usepackage{latexsym,amssymb}
\usepackage{comment}
\usepackage{enumerate}
\usepackage{graphicx}
\usepackage{authblk}
\usepackage{color}
\usepackage{hyperref}
\usepackage[french,english]{babel}
\usepackage[matrix,arrow]{xy}
\def \be{\begin{eqnarray*}}
\def \ee{\end{eqnarray*}}
\def \ben{\begin{enumerate}}
\def \een{\end{enumerate}}
\def \beit{\begin{itemize}}
\def \eeit{\end{itemize}}
\def \bui#1#2{\mathrel{\mathop{\kern 0pt#1}\limits^{#2}}}
\def \buil#1#2{\mathrel{\mathop{\kern 0pt#1}\limits_{#2}}}
\def \bfll{\begin{flushleft}}
\def \efll{\end{flushleft}}
\def \bflr{\begin{flushright}}
\def \eflr{\end{flushright}}

\def \findemo{\hfill$\square$\\}

\def \lra{\longrightarrow}

\def \R{\mathbb{R}}
\def \tr{\mathrm{tr}}

\def \Ric{\mathrm{Ric}}
\def \ric{\mathrm{ric}}

\newcommand{\rquot}[2]{\raisebox{0.5ex}{$#1$}\!/\!\raisebox{-0.5ex}{$#2$}}

\newtheorem{ethm}{Theorem}[section]

\newtheorem{elemme}[ethm]{Lemma}
\newtheorem{erem}[ethm]{Note}

\newtheorem{prop}[ethm]{Proposition}
\newtheorem{eexemple}[ethm]{Example}

\title{A generalised Ricci-Hessian
equation on Riemannian manifolds}

\author[1]{Nicolas Ginoux\thanks{\texttt{nicolas.ginoux@univ-lorraine.fr}}}
\author[1,2]{Georges Habib\thanks{\texttt{ghabib@ul.edu.lb}}}

\affil[1]{\footnotesize Universit\'e de Lorraine, CNRS, IECL, F-57000 Metz}
\affil[2]{\footnotesize Lebanese University, Faculty of Sciences II, Department 
of Mathematics, P.O. Box 90656 Fanar-Matn, Lebanon}

\begin{document}
\maketitle

\noindent\begin{center}\begin{tabular}{p{115mm}}
\begin{small}{\bf Abstract.}   In this paper, we prove new rigidity results related to some generalised Ricci-Hessian
equation on Riemannian manifolds.
\end{small}\\
\end{tabular}\end{center}

$ $\\

\noindent\begin{small}{\it Mathematics Subject Classification} (2010):  53C20, 58J60\\
\end{small}

\noindent\begin{small}{\it Keywords}: Obata equation, Tashiro equation, 
virtual Einstein equation\color{black}
\end{small}

\section{Introduction}\label{s:intro}

In this article, which follows \cite{GinouxHabibKathRicciHessian2018}, we continue investigating those Riemannian ma\-ni\-folds $(M^n,g)$ supporting a non-identically-vanishing function $f$ satisfying what we call the \emph{generalised Ricci-Hessian equation} \cite[Eq. (1)]{GinouxHabibKathRicciHessian2018}
\begin{equation}\label{eq:nabladf=-fric}
\nabla^2f=-f\cdot\Ric
\end{equation}
on $M$, where $\nabla^2f:=\nabla\nabla f$ denotes the Hessian of $f$ and $\Ric$ the Ricci-tensor of $(M^n,g)$, both seen as $(1,1)$-tensor fields.
Recall that this equation was first con\-si\-de\-red when studying the so-called \emph{skew-Killing-spinor-equation} \cite{GinouxHabibKathSKS2018}.
In \cite{GinouxHabibKathRicciHessian2018}, we proved that, provided sufficiently many symmetries preserving a solution $f$ are available on the underlying manifold $(M^n,g)$, only one of the following can occur: unless $f$ is constant and then $(M^n,g)$ is Ricci-flat, either $(M^n,g)$ is isometric to the Riemannian product of a real interval with a Ricci-flat manifold and $f$ is an affine-linear function on the interval; or $(M^n,g)$ is isometric to the Riemannian product of a Ricci-flat manifold with either the $2$-sphere or the hyperbolic plane and $f$ is the trivial extension of a solution to the Obata resp. Tashiro equation on the second factor.\\

In this article, we mainly show that, in many further situations, some of which are more general than those from \cite{GinouxHabibKathRicciHessian2018}, only those two possibilities can occur.\\

The article is structured as follows.
After preliminary remarks in Section \ref{s:prelim}, we describe and partially classify those warped products carrying solutions to (\ref{eq:nabladf=-fric}).
In Section \ref{s:dimWgeq2}, we turn to the case where the space of solutions to (\ref{eq:nabladf=-fric}) is at least $2$-dimensional.
Section \ref{s:homogcase} is dedicated to the homogeneous case, which remains partially open.
We conclude in Section \ref{s:Kaehlercase} with the case where the manifold is K\"ahler.\\

We underline that no full classification is available yet.
This is the object of future work.\\


{\bf Acknowledgment:}
Part of this work was done while the se\-cond-named author received the support 
of 
the Humboldt Foundation, the French embassy in Beirut via the SAFAR programme and the Alfried Krupp Wissenschaftskolleg, which he would like to thank.
We thank the Centre International de Rencontres Math\'ematiques (CIRM) where the article was finished.
We also thank Ines Kath for her support and interest in the first part of that work.

\color{black}
\section{Preliminary remarks}\label{s:prelim}


We start with preliminary results, some of which are already contained in \cite{GinouxHabibKathRicciHessian2018} but, for the sake of self-containedness, we give and reprove them here. \color{black}
From now on, we shall denote by $S$ the scalar curvature of $M$ and, for any 
function $h$ on $M$, by $\nabla h$ the gradient vector field of $h$ w.r.t. $g$ 
on $M$.
First observe that the equation $\nabla^2 f=-f\cdot\Ric$ is of course linear in 
$f$ but is also invariant under metric re\-sca\-ling: if $\overline{g}=\lambda^2 g$ 
for some nonzero real number $\lambda$, then $\overline{\nabla}^2 
f=\lambda^{-2}\overline{\nabla}^2 f$ (this comes from the re\-sca\-ling of 
the gradient) and $\overline{\Ric}=\lambda^{-2}\Ric$.
Let us denote by
\[W(M^n,g):=\left\{f\in C^\infty(M,\R)\,|\,\nabla^2f=-f\cdot\Ric\right\}\]
the real vector space of all smooth functions satisfying 
(\ref{eq:nabladf=-fric}) on $(M^n,g)$.

\begin{elemme}\label{l:eqnabladf=-fric}
Let $(M^n,g)$ be any connected Riemannian manifold carrying a smooth real-valued 
function $f$ satisfying {\rm(\ref{eq:nabladf=-fric})} on $M$.
\begin{enumerate}
\item\label{claim:Ricgradf} The gradient vector field $\nabla f$ of $f$ w.r.t. 
$g$ satisfies
\begin{equation}\label{eq:Ricnablaf}
\Ric(\nabla f)=\frac{S}{2}\nabla f+\frac{f}{4}\nabla S.
\end{equation}
\item\label{claim:fLaplacef} There exists a real constant $\mu$ such that 
\begin{equation}\label{eq:fLaplacef}
f\Delta f+2|\nabla f|^2=\mu.
\end{equation}

\item\label{claim:normRicci} The identity
\begin{equation}\label{eq:normRicci}
f|\Ric|^2=\frac{fS^2}{2}-\frac{1}{4}\langle\nabla f,\nabla 
S\rangle+\frac{f}{4}\Delta S
\end{equation}
holds on $M$.
\item\label{claim:confchange} If $n>2$ and $f$ is everywhere positive or 
negative, then $f$ solves 
{\rm(\ref{eq:nabladf=-fric})} if and only if, setting $u:=\frac{1}{2-n}\ln|f|$, 
the metric $\overline{g}:=e^{2u}g$ satisfies 
$\overline{\ric}=(\overline{\Delta}u)\overline{g}-(n-2)(n-3)du\otimes du$ on 
$M$ 
and in that case $\overline{\Delta}u=-\frac{\mu}{n-2}e^{2(n-3)u}$.
In particular, if $n=3$, the existence of a positive solution $f$ to 
{\rm(\ref{eq:nabladf=-fric})} is equivalent to $(M,f^{-2}g)$ being Einstein 
with 
scalar curvature $-3\overline{\Delta}\ln|f|$.
\item\label{claim:Mclosedfconstantsign} If $M$ is closed and $f$ is everywhere 
positive or negative, then $f$ is 
constant on $M$.
\item\label{claim:N0totgeod}
If nonempty, the vanishing set 
$N_0:=f^{-1}(\{0\})$ of $f$ is a 
scalar-flat totally geodesic hypersurface of $M$.
\item\label{claim:RXYnablaf}  For any $x\in M$ and all $X,Y\in
T_xM$, the identity 
\begin{equation}\label{eq:RXYnablaf}
R_{X,Y}\nabla 
f=-X(f)\mathrm{Ric}(Y)+Y(f)\mathrm{Ric}(X)-f\left((\nabla_X\mathrm{Ric} 
)Y-(\nabla_Y\mathrm{Ric})X\right)
\end{equation}
holds on $M$.
As a consequence, at any critical point of $f$, the Ricci-tensor must be 
Codazzi.
\item\label{claim:dimWleqn+1} The dimension of $W(M^n,g)$ is always
at most $n+1$.
\item\label{claim:MEinstein} If furthermore $M$ is Einstein or
$2$-dimensional, then $M$ is Ricci-flat or
$n=2$ and in that case $M$ has constant curvature.
In particular, when $(M^2,g)$ is complete, there exists a nonconstant function 
$f$ satisfying {\rm(\ref{eq:nabladf=-fric})} if and only if, up to re\-sca\-ling 
the 
metric, the manifold $(M^2,g)$ is isometric to either the round sphere 
$\mathbb{S}^2$ and $f$ is a nonzero eigenfunction associated to the first 
positive Laplace eigenvalue; or to flat $\R^2$ or cylinder 
$\mathbb{S}^1\times\R$ and $f$ is an affine-linear function; or to the 
hyperbolic plane $\mathbb{H}^2$ and $f$ is a solution to the Tashiro equation 
$\nabla^2f=f\cdot\mathrm{Id}$.
\item\label{claim:Sconstantnugeodesic} If $S$ is constant, then outside the set 
of critical points of $f$, the 
vector field $\nu:=\frac{\nabla f}{|\nabla f|}$ is geodesic.
Moreover, assuming $(M^n,g)$ to be also complete,
\begin{enumerate}
\item\label{claim:SconstantSpositivecritf} if $S>0$, then up to re\-sca\-ling the 
metric as well as $f$, we may assume 
that $S=2$ and that $\mu=f\Delta f+2|\nabla f|^2=2$ on $M$, in which case the 
function $f$ has $1$ as maximum and $-1$ as minimum value and those are the 
only 
critical values of $f$;
\item\label{claim:SconstantSequal0}if $S=0$  and $f$ is 
nonconstant, then $(M^n,g)$ is Ricci-flat, in 
particular it is isometric to 
$(\R\times\Sigma,dt^2\oplus g_{\Sigma})$ for some complete Ricci-flat 
Riemannian 
manifold $(\Sigma,g_{\Sigma})$ and, up to reparametrization, the function $f$ 
is 
given by $f(t,x)=t$ for all $(t,x)\in\R\times\Sigma$;
\item\label{claim:SconstantSnegativecritf} if $S<0$, then up to re\-sca\-ling the 
metric, we may assume that $S=-2$ on 
$M$, in which case one of the following holds:
\begin{enumerate}
\item\label{claim:SconstantSnegativemupositivecritf} if $\mu>0$, then up to 
re\-sca\-ling $f$ we may assume that $\mu=2$, in which 
case $f$ has no critical value and $f(M)=\R$, in particular $M$ is noncompact;
\item\label{claim:SconstantSnegativemuzerocritf} if $\mu=0$, then $f$ has no 
critical value and empty vanishing set and, up 
to changing $f$ into $-f$, we have $f(M)=(0,\infty)$, in particular $M$ is 
noncompact;
\item\label{claim:SconstantSnegativemunegativecritf} if $\mu<0$, then up to 
re\-sca\-ling $f$ we may assume that $\mu=-2$, in which 
case $f$ has a unique critical value, which, up to changing $f$ into $-f$, can 
be assumed to be a minimum; moreover, $f(M)=[1,\infty)$, in particular $M$ is 
noncompact.
\end{enumerate}

\end{enumerate}
\end{enumerate}
\end{elemme}


{\it Proof:}
The proof of statement \ref{claim:Ricgradf}. follows that of \cite[Lemma 
4]{KimKim03}.
On the one hand, we take the codifferential of $\nabla^2f$ and obtain, choosing 
a local orthonormal basis $(e_j)_{1\leq j\leq n}$ of $TM$ and using the 
Weitzenb\"ock 
formula for $1$-forms:
\begin{eqnarray}
\nonumber\delta\nabla^2f&=&-\sum_{j=1}^n\left(\nabla_{e_j}\nabla^2f\right)(e_j)\\
\nonumber&=&-\sum_{j=1}^n\left(\nabla_{e_j}\nabla_{e_j}\nabla f-\nabla_{\nabla_{e_j}e_j}\nabla
f\right)\\
\nonumber&=&\nabla^*\nabla (\nabla f)\\
\label{eq:Weitzenboeck}&=&\Delta(\nabla f)-\Ric(\nabla f).
\end{eqnarray}
On the other hand, by (\ref{eq:nabladf=-fric}) and the formula 
$\delta\Ric=-\frac{1}{2}\nabla S$,
\begin{eqnarray*}
\delta\nabla^2f&=&\delta\left(-f\cdot\Ric\right)\\
&=&\Ric(\nabla f)-f\cdot\delta\Ric\\
&=&\Ric(\nabla f)+\frac{f}{2}\nabla S.
\end{eqnarray*}
Comparing both identities, we deduce that $\Delta(\nabla f)=2\Ric(\nabla 
f)+\frac{f}{2}\nabla S$.
But identity (\ref{eq:nabladf=-fric}) also gives 
\begin{equation}\label{eq:DeltaffS}\Delta 
f=-\mathrm{tr}\left(\nabla^2f\right)=fS,\end{equation}
so that $\Delta(\nabla f)=\nabla (\Delta f)=\nabla(f S)=S\nabla f+f\nabla S$ 
and 
therefore $\Ric(\nabla f)=\frac{S}{2}\nabla f+\frac{f}{4}\nabla S$, which is 
(\ref{eq:Ricnablaf}).\\
By (\ref{eq:nabladf=-fric}) and (\ref{eq:Ricnablaf}), we have
\begin{eqnarray*}
2\nabla(|\nabla f|^2)&=&4\nabla_{\nabla f}^2f\\
&=&-4f\cdot\Ric(\nabla f)\\
&=&-4f\cdot\left(\frac{S}{2}\nabla f+\frac{f}{4}\nabla S\right)\\
&=&-2Sf\nabla f-f^2\nabla S\\
&=&-\nabla(Sf^2)\\
&\bui{=}{\rm(\ref{eq:DeltaffS})}&-\nabla(f\Delta f),
\end{eqnarray*}
from which (\ref{eq:fLaplacef}) follows.\\

Taking the codifferential of (\ref{eq:Ricnablaf}), we obtain on the one 
hand, using $\delta\Ric=-\frac{1}{2}\nabla S$:
\begin{eqnarray*}
\delta(\Ric\nabla f)&=&\langle\delta\Ric,\nabla f\rangle-\langle\Ric,\nabla^2 
f\rangle\\
&\bui{=}{\rm(\ref{eq:nabladf=-fric})}&-\frac{1}{2}\langle\nabla S,\nabla 
f\rangle+f|\Ric|^2.
\end{eqnarray*}
On the other hand, the codifferential of the r.h.s. of (\ref{eq:Ricnablaf}) is 
given by
\begin{eqnarray*}
\delta(\frac{S}{2}\nabla f+\frac{f}{4}\nabla S)&=&-\frac{1}{2}\langle\nabla 
 S,\nabla f\rangle+\frac{S}{2}\Delta f-\frac{1}{4}\langle\nabla f,\nabla 
S\rangle+\frac{f}{4}\Delta S\\
&=&-\frac{3}{4}\langle\nabla f,\nabla S\rangle+\frac{S^2f}{2}+\frac{f}{4}\Delta 
S.
\end{eqnarray*}
Comparing both identities yields (\ref{eq:normRicci}).\\

If $f$ vanishes nowhere, then up to changing $f$ into $-f$, we may assume that 
$f>0$ on $M$.
Writing $f$ as $e^{(2-n)u}$ for some real-valued function $u$ (that is, 
$u=\frac{1}{2-n}\ln f$), the Ricci-curvatures (as $(0,2
)$-tensor fields) $\ric$ and $\overline{\ric}$ of $(M,g)$ and 
$(M,\overline{g}=e^{2u}g)$ respectively are related as follows:
\begin{equation}\label{eq:Ricciggbar}
\overline{\ric}=\ric+(2-n)(\nabla du-du\otimes du)+(\Delta u-(n-2)|du|_g^2)g.
\end{equation}
But $\nabla df=(n-2)^2f\cdot du\otimes du+(2-n)f\cdot\nabla du$ and the Laplace 
operators $\Delta$ of $(M,g)$ and $\overline{\Delta}$ of $(M,\overline{g})$ are 
related via $\overline{\Delta}v=e^{-2u}\cdot(\Delta v-(n-2)g(du,dv))$ for any 
function $v$, so that 
\begin{eqnarray*}
\overline{\ric}&=&\ric+\frac{1}{f}\nabla df-(n-2)^2du\otimes du+(n-2) du\otimes 
du+(\overline{\Delta}u)\overline{g}\\
&=&\ric+\frac{1}{f}\nabla df-(n-2)(n-3)du\otimes 
du+(\overline{\Delta}u)\overline{g}.
\end{eqnarray*}
As a consequence, $f$ satisfies (\ref{eq:nabladf=-fric}) if and only if 
$\overline{\ric}=(\overline{\Delta}u)\overline{g}-(n-2)(n-3)du\otimes du$ holds 
on $M$.
Moreover, 
\begin{eqnarray*}
f\Delta f+2|df|_g^2&=&f\cdot\left(-(n-2)^2f|du|_g^2-(n-2)f\Delta 
u\right)+2(n-2)^2f^2|du|_g^2\\
&=&-(n-2)f^2\cdot\left(\Delta u-(n-2)|du|_g^2\right)\\
&=&-(n-2)f^2\cdot e^{2u}\cdot\overline{\Delta}u\\
&=&-(n-2)e^{2(2-n)u}\cdot e^{2u}\cdot\overline{\Delta}u\\
&=&-(n-2)e^{2(3-n)u}\cdot\overline{\Delta}u,
\end{eqnarray*}
in particular (\ref{eq:fLaplacef}) yields 
$\overline{\Delta}u=-\frac{\mu}{n-2}e^{2(n-3)u}$.
In dimension $3$, we notice that $\overline{\Delta}u=\frac{\overline{S}}{3}$.
This shows statement \ref{claim:confchange}.\\
If $f$ vanishes nowhere, then again we may assume that $f>0$ on $M$.
Since $M$ is closed, $f$ has a minimum and a maximum.
At a point $x$ where $f$ attains its maximum, we have $\mu=f(x)(\Delta 
f)(x)+2|\nabla_xf|^2=f(x)(\Delta f)(x)\geq0$.
In the same way, $\mu=f(y)(\Delta f)(y)\leq0$ at any point $y$ where $f$ 
attains 
its minimum.
We deduce that $\mu=0$ which, by integrating the identity $f\Delta f+2|\nabla 
f|^2=\mu$ on $M$, yields $df=0$.
This shows statement \ref{claim:Mclosedfconstantsign}.\\
The first part of statement \ref{claim:N0totgeod}. is the consequence of the 
following very 
general fact \cite[Prop. 1.2]{HePetersenWylie11102455}, that we state and 
reprove here for the sake of completeness: if some smooth real-valued function 
$f$ satisfies $\nabla^2f=fq$ for some quadratic form $q$ on $M$, then the 
subset 
$N_0=f^{-1}\left(\{0\}\right)$ is -- if nonempty -- a totally geodesic smooth 
hypersurface of $M$.
First, it is a smooth hypersurface because of $d_xf\neq0$ for all $x\in N_0$: 
namely if $c\colon\mathbb{R}\to M$ is any geodesic with $c(0)=x$, then the 
function $y:=f\circ c$ satisfies the se\-cond order linear ODE 
$y''=\langle\nabla_{\dot{c}}^2f,\dot{c}\rangle=q(\dot{c},\dot{c})\cdot y$ on 
$\mathbb{R}$ with the initial condition $y(0)=0$; if $d_xf=0$, then $y'(0)=0$ 
and hence $y=0$ on $\mathbb{R}$, which would imply that $f=0$ on $M$ by 
geodesic 
connectedness, contradiction.
To compute the shape operator $W$ of $N_0$ in $M$, we define $\nu:=\frac{\nabla 
f}{|\nabla f|}$ to be a unit normal to $N_0$.
Then for all $x\in N_0$ and $X\in T_xM$,
\begin{eqnarray}\label{eq:nablaXnu}
\nonumber\nabla_X^M\mathrm{\nu}&=&X\left(\frac{1}{|\nabla f|}\right)\cdot\nabla 
f+\frac{1}{|\nabla f|}\cdot\nabla_X^M\nabla f\\
\nonumber&=&-\frac{X\left(|\nabla f|^2\right)}{2|\nabla f|^3}\cdot\nabla 
f+\frac{1}{|\nabla f|}\cdot\nabla_X^M\nabla f\\
&=&\frac{1}{|\nabla f|}\cdot\left(\nabla_X^2 f-\langle\nabla_X^2 
f,\nu\rangle\cdot\nu\right),
\end{eqnarray}
in particular $W_x=-(\nabla\nu)_x=0$ because of 
$\left(\nabla^2f\right)_x=f(x)q_x=0$.
This shows that $N_0$ lies totally geodesically in $M$.\\
Now recall Gau\ss{} equations for Ricci curvature: for every $X\in TN_0$,
$$\Ric_{N_0}(X)=\Ric(X)^T-R_{X,\nu}^M\nu+\mathrm{tr}_g(W)\cdot WX-W^2X,$$
where $\Ric(X)^T=\Ric(X)-\ric(X,\nu)\nu$ is the component of the Ricci 
curvature 
that is tangential to the hypersurface $N_0$.
As a straightforward consequence, if $S_{N_0}$ denotes the scalar curvature of 
$N_0$,
$$S_{N_0}=S-2\ric(\nu,\nu)+\left(\mathrm{tr}_g(W)\right)^2-|W|^2.$$
Here, $W=0$ and $\Ric(\nu)=\frac{S}{2}\nu$ along $N_0$ because $N_0$ lies 
totally geodesically in $M$, so that
$$S_{N_0}=S-2\ric(\nu,\nu)=S-S=0.$$
This proves $N_0$ to be scalar-flat and statement \ref{claim:N0totgeod}.\\
 As for claim \ref{claim:RXYnablaf}., a straightforward
consequence 
of (\ref{eq:nabladf=-fric}) is that, at every $x\in M$ 
and for all $X,Y\in T_xM$, we have
\be 
R_{X,Y}\nabla f&=&\left[\nabla_X,\nabla_Y\right]\nabla f-\nabla_{[X,Y]}^2f\\
&=&-X(f)\mathrm{Ric}(Y)+Y(f)\mathrm{Ric}(X)-f\left((\nabla_X\mathrm{Ric}
)Y-(\nabla_Y\mathrm{Ric})X\right),
\ee
which is identity (\ref{eq:RXYnablaf}).
In particular, because $0$ cannot be a critical value of $f$ by statement 
\ref{claim:N0totgeod}., the Ricci-tensor of $(M^n,g)$ must be Codazzi at every 
critical point of $f$.
This proves claim \ref{claim:RXYnablaf}.\\
Statement \ref{claim:dimWleqn+1}., which can be found in \cite[Prop. 
1.1]{HePetersenWylie11102455}, is a further consequence of the 
general fact mentioned above that any $f\in W(M^n,g)$ is uniquely determined by 
its 
value as well as the value of its gradient at a given point.
This implies that, given any $x\in M$, the linear map
\be 
W(M^n,g)&\longrightarrow&\R\times T_xM\\
f&\longmapsto&(f(x),(\nabla f)(x))
\ee
is injective, which proves claim \ref{claim:dimWleqn+1}.
Note that the upper bound $n+1$ for $\dim(W(M^n,g))$ is obviously attained when 
$(M^n,g)=(\R^n,\mathrm{can})$ is the flat Euclidean space.\\
Statement \ref{claim:MEinstein}. can be considered as standard.
In dimension $2$, $\Ric=\frac{S}{2}\mathrm{Id}=K\mathrm{Id}$, 
where 
$K$ is the Gau\ss{} curvature of $(M^2,g)$.
But we also know that $\Ric(\nabla f)=\frac{S}{2}\nabla f+\frac{f}{4}\nabla 
S=K\nabla f+\frac{f}{2}\nabla K$.
Comparing both identities and using the fact that $\{f\neq0\}$ is dense in $M$ 
leads to $\nabla K=0$, that is, $M$ has constant Gau\ss{} curvature.
Up to re\-sca\-ling the metric as well as $f$, we may assume that 
$S,\mu\in\{-2,0,2\}$.
If $M^2$ is complete with constant $S>0$ (hence $K=1$) and $f$ is nonconstant, 
then $\mu>0$ so that, by Obata's solution to the equation $\nabla^2 
f+f\cdot\mathrm{Id}_{TM}=0$, the manifold $M$ must be isometric to the round 
sphere of radius $1$ and the function $f$ must be a nonzero eigenfunction 
associated to the first positive eigenvalue of the Laplace operator on 
$\mathbb{S}^2$, see \cite[Theorem A]{Obata62}.
If $M^2$ is complete and has vanishing curvature, then its universal cover is 
the flat $\R^2$ and the lift $\tilde{f}$ of $f$ to $\R^2$ must be an 
affine-linear function of the form $\tilde{f}(x)=\langle a,x\rangle+b$ for some 
nonzero $a\in\R^2$ and some $b\in\R$; since the only possible nontrivial 
quotients of $\R^2$ on which $\tilde{f}$ may descend are of the form 
$\rquot{\R}{\mathbb{Z}\cdot \check{a}}\times\R$ for some nonzero $\check{a}\in 
a^\perp$, the manifold $M$ itself must be either flat $\R^2$ or such a flat 
cylinder.
If $M^2$ is complete with constant $S<0$, then $f$ satisfies the Tashiro 
equation $\nabla^2f=f\cdot\mathrm{Id}_{TM}$.
But then Y.~Tashiro proved that $(M^2,g)$ must be isometric to the hyperbolic 
plane of constant sectional curvature $-1$, see e.g. \cite[Theorem 2 
p.252]{Tashiro65}, see also \cite[Theorem G]{Kanai1983}.
Note that the functions $f$ listed above on $\mathbb{S}^2$, $\mathbb{R}^2$, 
$\mathbb{S}^1\times\mathbb{R}$ or $\mathbb{H}^2$ obviously satisfy 
(\ref{eq:nabladf=-fric}).\\
If $(M^n,g)$ is Einstein with $n\geq3$, then it has constant scalar curvature 
$S$
and $\Ric=\frac{S}{n}\cdot\mathrm{Id}$.
But again the identity $\Ric(\nabla f)=\frac{S}{2}\nabla f+\frac{f}{4}\nabla 
S=\frac{S}{2}\nabla f$ yields $n=2$ unless $S=0$ and thus $M$ is Ricci-flat.
Therefore, $n=2$ is the only possibility for non-Ricci-flat Einstein $M$.
This shows statement \ref{claim:MEinstein}.\\
If $S$ is constant, then $\Ric(\nabla f)=\frac{S}{2} \nabla f$.
As a consequence, $\nabla_{\nabla f}^2f=-f\Ric(\nabla f)=-\frac{Sf}{2}\nabla f$.
But, as already observed in e.g. \cite[Prop. 1]{RanjSanth95}, away from its 
vanishing set, the gradient of $f$ is a pointwise eigenvector of its Hessian if 
and only if the vector field $\nu=\frac{\nabla f}{|\nabla f|}$ is geodesic, see 
(\ref{eq:nablaXnu}) above.
Assuming furthermore $(M^n,g)$ to be complete, we can rescale as before $f$ and 
$g$ such that $S,\mu\in\{-2,0,2\}$.
In case $S>0$ and hence $S=2$, necessarily $\mu>0$ holds and thus $\mu=2$.
But then $f^2+|\nabla f|^2=1$, so that the only critical points of $f$ are 
those 
where $f^2=1$, which by $f^2\leq1$ shows that the only critical points of $f$ 
are those where $f=\pm1$ and hence where $f$ takes a maximum or minimum value.
Outside critical points of $f$, we may consider the function $y:=f\circ 
\gamma\colon \R\to\R$, where $\gamma\colon \R\to M$ is a maximal integral curve 
of the geodesic vector field $\nu$.
Then $y$ satisfies $y'=|\nabla f|\circ\gamma>0$ and $y(t)^2+y'(t)^2=1$, so that 
$y'=\sqrt{1-y^2}$ and therefore there exists some $\phi\in\R$ such that 
$$y(t)=\cos(t+\phi)\qquad\forall\,t\in \R.$$
Since that function obviously changes sign and $0$ is not a critical value of 
$f$, we can already deduce that $f$ changes sign, in particular 
$N_0=f^{-1}(\{0\})$ is nonempty.
Moreover, the explicit formula for $y$ shows that $f$ must have critical 
points, 
which are precisely those where $\cos$ reaches its minimum or maximum value.
This shows statement 
\ref{claim:SconstantSpositivecritf}.\\

In case $S=0$, we have $\Ric=0$ by (\ref{eq:normRicci}) since $f$ is assumed to 
be nonconstant.
This proves statement \ref{claim:SconstantSequal0}.\\

In case $S<0$ and thus $S=-2$, there are still three possibilities for $\mu$:
\begin{enumerate}[$\bullet$]
\item If $\mu>0$, then $\mu=2$ and (\ref{eq:fLaplacef}) becomes $-f^2+|\nabla 
f|^2=1$, hence $f$ has no critical point.
If $\gamma$ is any integral curve of the normalised gradient vector field
$\nu=\frac{\nabla f}{|\nabla f|}$, then the function $y:=f\circ\gamma$ 
satisfies 
the ODEs $y'=\sqrt{1+y^2}$, therefore $y(t)=\sinh(t+\phi)$ for some real 
constant $\phi$. 
In particular, $f(M)=\R$ and $M$ cannot be compact.
\item If $\mu=0$, then (\ref{eq:fLaplacef}) becomes $f^2=|\nabla f|^2$.
But since no point where $f$ vanishes can be a critical point by the fifth 
statement, $f$ has no critical point and therefore must be of constant sign.
Up to turning $f$ into $-f$, we may assume that $f>0$ and thus $f=|\nabla f|$.
Along any integral curve $\gamma$ of $\nu=\frac{\nabla f}{|\nabla f|}$, the 
function $y:=f\circ\gamma$ satisfies $y'=y$ and hence $y(t)=C\cdot e^t$ for 
some 
positive constant $C$.
This shows $f(M)=(0,\infty)$, in particular $M$ cannot be compact.
\item If $\mu<0$, then $\mu=-2$ and (\ref{eq:fLaplacef}) becomes $-f^2+|\nabla 
f|^2=-1$.
As a consequence, because of $f^2=1+|\nabla f|^2\geq1$, the function $f$ has 
constant sign and hence we may assume that $f\geq1$ up to changing $f$ into 
$-f$.
In particular, the only possible critical value of $f$ is $1$, which is an 
absolute minimum of $f$.
If $\gamma$ is any integral curve of the normalised gradient vector field
$\nu=\frac{\nabla f}{|\nabla f|}$, which is defined at least on the set of 
regular points of $f$, then the function $y:=f\circ\gamma$ satisfies the ODEs 
$y'=\sqrt{y^2-1}$, therefore $y(t)=\cosh(t+\phi)$ for some real constant 
$\phi$. 
Since that function has an absolute minimum, it must have a critical point.
It remains to notice that $f(M)=[1,\infty)$ and thus that $M$ cannot be compact.
\end{enumerate}
This shows statement \ref{claim:SconstantSnegativecritf}.
\findemo

\color{black}

\begin{eexemple}\label{ex:nabladf=-fricdim3}
{\rm In dimension $3$, Lemma \ref{l:eqnabladf=-fric} implies that, starting 
with 
any Einstein manifold
 -- or, equivalently, any manifold with constant sectional curvature -- $(M^3,g)$ and any real function $u$ such that $\Delta u=\frac{S}{3}$, the
function $f:=e^{-u}$ satisfies (\ref{eq:nabladf=-fric}) on the manifold 
$(M,\overline{g}=e^{-2u}g)$.
In particular, since there is an infinite-dimensional space of harmonic 
functions on any nonempty open subset $M$ of $\R^3$, there are many 
nonhomothetic conformal metrics on such $M$ for which nonconstant solutions of 
(\ref{eq:nabladf=-fric}) exist.
{ As a first consequence, there exist metrics with nonconstant scalar curvature 
on $\R^3$ for which there are nonconstant solutions of (\ref{eq:nabladf=-fric}).
However it remains unclear whether such metrics can be \emph{complete} or not.}
On any nonempty open subset of the $3$-dimensional hyperbolic space 
$\mathbb{H}^3$ with constant sectional curvature $-1$, there is also an 
infinite-dimensional affine space of solutions to the Poisson equation $\Delta 
u=-2${: in geodesic polar coordinates about any fixed point $p\in\mathbb{H}^3$, 
assuming $u$ to depend only on the geodesic distance $r$ to $p$, that Poisson 
equation is a se\-cond-order linear ODE in $u(r)$ and therefore has infinitely 
many affinely independent solutions}.
In particular, there are also lots of conformal metrics on $\mathbb{H}^3$ for 
which nonconstant solutions of (\ref{eq:nabladf=-fric}) exist.\\
{ Note however that, although $\mathbb{H}^3$ is conformally equivalent to the 
unit open ball $\mathbb{B}^3$ in $\R^3$, we do not obtain the same solutions to 
the equation depending on the metric we start from.
Namely, we can construct solutions of (\ref{eq:nabladf=-fric}) starting from 
the 
Euclidean metric $g$ and from the hyperbolic metric $e^{2w}g$ on 
$\mathbb{B}^3$, 
where $e^{2w(x)}=\frac{4}{(1-|x|^2)^2}$ at any $x\in\mathbb{B}^3$.
In both cases we obtain solutions of (\ref{eq:nabladf=-fric}) by conformal 
change of the metric.
Since $g$ and $e^{2w}g$ lie in the same conformal class, the question arises 
whether solutions coming from $e^{2w}g$ can coincide with solutions coming from 
$g$ on $\mathbb{B}^3$.
Assume $f$ were a solution of (\ref{eq:nabladf=-fric}) arising by conformal 
change of $g$ (by $e^{-2u}$ for some $u\in C^\infty(\mathbb{B}^3)$) and by 
conformal change of $e^{2w}g$ (by $e^{-2v}$ for some $v\in 
C^\infty(\mathbb{B}^3)$).
Then $f=e^{-u}=e^{-v}$ and thus $v=u$ would hold, therefore $u$ would satisfy 
$\Delta_g u=0$ as well as $\Delta_{e^{2w}g}u=-2$, where $\Delta_h$ is the Laplace operator associated to the metric $h$.
In particular
\begin{eqnarray*}
0&=&\Delta_g u\\
&=&e^{2w}\Delta_{e^{2w}g}u+\langle dw,du\rangle_g\\
&=&-2e^{2w}+\langle dw,du\rangle_g.
\end{eqnarray*}
But the r.h.s. of the last identity has no reason to vanish in general.
Note also that the conformal metrics themselves have no reason to coincide, 
since otherwise $e^{-2v}e^{2w}g=e^{-2u}g$ would hold hence $u=v-w$ as well and 
the same kind of argument would lead to an equation that is generally not 
fulfilled.
}
}
\end{eexemple}

\begin{erem}\label{r:critptf}
\noindent{\rm
{If $S$ is a nonzero constant and $M$ is closed, then the function $f$ is an 
eigenfunction for the scalar Laplace operator associated to the eigenvalue $S$ 
on $(M,g)$ and it has at least two nodal domains.
Mind however that $S$ is not necessarily the first positive Laplace eigenvalue 
on $(M,g)$.
E.g. consider the Riemannian manifold $M=\mathbb{S}^2\times\Sigma^{n-2}$ which 
is the product of standard $\mathbb{S}^2$ with a closed Ricci-flat manifold 
$\Sigma^{n-2}$, then the first positive Laplace eigenvalue of $\Sigma$ can be 
made arbitrarily small by re\-sca\-ling its metric; since the Laplace spectrum of 
$M$ is the sum of the Laplace spectra of $\mathbb{S}^2$ and $\Sigma$, the first 
Laplace eigenvalue on $M$ can be made as close to $0$ as desired by re\-sca\-ling 
the metric on $\Sigma$.
}
}
\end{erem}

Next we give a closer look at the case where the scalar curvature of $(M^n,g)$ is constant.

\begin{prop}\label{l:eqnabladf=-fricscalconstant}
Let $(M^n,g)$ be any connected Riemannian manifold carrying a nonzero smooth real-valued
function $f$ satisfying {\rm(\ref{eq:nabladf=-fric})} on $M$.
Assume the scalar curvature $S$ of $(M^n,g)$ to be constant and nonvanishing.
Up to re\-sca\-ling the metric $g$ on
$M$
it may be assumed  that $S=2\varepsilon$ for
some $\varepsilon\in\{\pm1\}$.\\
Then the following holds.
\begin{enumerate}
\item\label{statement:Ncscal0} Every regular level hypersurface
$N_c:=f^{-1}(\{c\})$ of $f$ must have vanishing scalar curvature and its
Ricci-tensor be given by $\Ric_{N_c}=-\frac{f}{|\nabla
f|^2}\left(\nabla_{\nabla
f}\Ric\right)$.
\item\label{statement:valpRicci} If either $n=3$ or both $n\geq4$ and $\Ric$ is assumed to
be nonnegative when $\varepsilon=1$ resp. nonpositive when $\varepsilon=-1$,
then the Ricci-tensor has pointwise $2$ eigenvalues, $\varepsilon$ with
multiplicity $2$ and $0$ with multiplicity $n-2$.
\item\label{statement:scalcstn=3} If $n=3$, the manifold $(M^3,g)$ must be
isometric to either
$S^2(\varepsilon)\times\R$ or $S^2(\varepsilon)\times S^1$ with product
Riemannian metric, where $S^2(\varepsilon)$ is the simply-connected complete
surface of constant curvature
$\varepsilon\in\{\pm1\}$; and $f$ must be the trivial extension to $M$ of a
solution
of the Obata resp.
Tashiro equation on $S^2(1)=\mathbb{S}^2$ (if $\varepsilon=1$) resp. $S^2(-1)=\mathbb{H}^2$
(if
$\varepsilon=-1$).
\end{enumerate}
\end{prop}

{\it Proof:}
We look at the Gau\ss{} equations for Ricci and scalar curvature  along each
$N_c:=f^{-1}(\{c\})$ for any regular value $c$ of $f$.
Denoting $W=-\nabla \nu=\frac{f}{|\nabla f|}\Ric^T=\frac{f}{|\nabla f|}\Ric$
the
Weingarten-endomorphism-field of $N_c$ in $M$, where $\Ric^T$ is the pointwise
orthogonal projection of $\Ric$ onto $TN_c$, we have
$\mathrm{tr}(W)=\frac{f}{|\nabla f|}\cdot\frac{S}{2}$ by
$\Ric(\nu)=\frac{S}{2}\nu$.
As a consequence, we have, for all $X\in TN_c$:
\begin{eqnarray*}
\Ric(X)&=&\Ric(X)^T\\
&=&\Ric_{N_c}(X)+W^2X-\mathrm{tr}(W)WX+R_{X,\nu}\nu\\
&=&\Ric_{N_c}(X)+\frac{f^2}{|\nabla
f|^2}\left(\Ric^2(X)-\frac{S}{2}\Ric(X)\right)+R_{X,\nu}\nu.
\end{eqnarray*}
But we can compute the curvature term $R_{X,\nu}\nu$ explicitly from
(\ref{eq:RXYnablaf}): for any $X\in TN_c$,
\begin{eqnarray}\label{eq:RXnunu}
\nonumber R_{X,\nu}\nu&=&-\frac{X(f)}{|\nabla
f|}\mathrm{Ric}(Y)+\frac{\nu(f)}{|\nabla f|}\mathrm{Ric}
(X)-\frac{f}{|\nabla f|}\left((\nabla_X\mathrm { Ric}
)\nu-(\nabla_\nu\mathrm{Ric})X\right)\\
\nonumber &=&\Ric(X)-\frac{f}{|\nabla
f|}\left(\nabla_X(\underbrace{\Ric\nu}_{\frac{S}{
2}\nu})-\Ric(\nabla_X\nu)\right)+\frac { f } { |\nabla
f|}(\nabla_\nu\mathrm{Ric})X\\
\nonumber &=&\Ric(X)+\frac{f}{|\nabla
f|}\left(\frac{S}{2}\mathrm{Id}-\Ric\right)(WX)+\frac{f}{|\nabla
f|}(\nabla_\nu\mathrm{Ric})X\\
&=&\Ric(X)+\frac{f^2}{|\nabla
f|^2}\left(\frac{S}{2}\Ric(X)-\Ric^2(X)\right)+\frac{f}{|\nabla
f|}(\nabla_\nu\mathrm{Ric})X,
\end{eqnarray}
so that, with
$(\nabla_\nu\Ric)(\nu)=\nabla_\nu(\Ric(\nu))-\Ric(\nabla_\nu\nu)=\nabla_\nu(
\frac{S}{2}\nu)=0$ on $M$, we obtain
\[\Ric_{N_c}=-\frac{f}{|\nabla
f|}\cdot\nabla_\nu\mathrm{Ric},\]
as claimed in statement \ref{statement:Ncscal0}.
That identity has important consequences.
First, choosing a local o.n.b. $(e_j)_{1\leq j\leq n-1}$ of $TN_c$,
\begin{eqnarray*}
S_{N_c}&=&\sum_{j=1}^{n-1}\langle\Ric_{N_c}(e_j),e_j\rangle\\
&=&-\frac{f}{|\nabla
f|}\cdot\sum_{j=1}^{n-1}\langle(\nabla_\nu\Ric)(e_j),e_j\rangle\\
&=&-\frac{f}{|\nabla
f|}\cdot\left(\sum_{j=1}^{n-1}\langle(\nabla_\nu\Ric)(e_j),
e_j\rangle+\langle(\nabla_\nu\Ric)(\nu),\nu\rangle\right)\\
&&+\frac{f}{|\nabla
f|}\cdot \langle(\nabla_\nu\Ric)(\nu),\nu\rangle,
\end{eqnarray*}
so that
$$S_{N_c}=-\frac{f}{|\nabla
f|}\cdot\mathrm{tr}(\nabla_\nu\Ric)=-\frac{f}{|\nabla
f|}\cdot\nu(\mathrm{tr}(\Ric))=-\frac{f}{|\nabla f|}\cdot\nu(S)=0.$$
Therefore, each level hypersurface $N_c$ is scalar-flat.
This concludes the proof of statement \ref{statement:Ncscal0}.
We turn to \ref{statement:valpRicci}.
Because of $S$ being constant, we already know by (\ref{eq:Ricnablaf}) that,
outside its vanishing set, the gradient vector field $\nabla f$ of $f$ is a
pointwise eigenvector for the Ricci tensor associated to the eigenvalue
$\frac{S}{2}=\varepsilon$.
Writing the Ricci tensor as $\Ric=\varepsilon\nu^\flat\otimes\nu+\Ric^T$, where
$\Ric^T$ is a pointwise symmetric endomorphism of $\nu^\perp\subset TM$, we
deduce from (\ref{eq:normRicci}) and the fact that $\{f\neq0\}$ is dense in $M$
that
\begin{equation}\label{eq:normRicTscalconstant}
|\Ric^T|^2=\frac{S^2}{4}=1
\end{equation}
on $\{\nabla f\neq0\}$.
Since $\tr(\Ric^T)=\frac{S}{2}=\varepsilon$, identity
(\ref{eq:normRicTscalconstant}) implies that, outside the critical set, the set
of possible pointwise
eigenvalues of $\Ric^T$ {\sl a priori} stands in one-to-one correspondence with
the sphere $\mathbb{S}^{n-3}$ of
dimension $n-3$.
If $n=3$, then this means that $\Ric^T$ has pointwise the
eigenvalues $\varepsilon$ and $0$, each of multiplicity one, on the regular set of
$f$.
If $n\geq4$, we assume furthermore that $\Ric\geq0$ when $\varepsilon=1$ and $\Ric\leq0$
when
$\varepsilon=-1$.
In that case, (\ref{eq:normRicTscalconstant}) implies that $\Ric^T$ has
exactly one eigenvalue that is equal to $\varepsilon$ and that all other
eigenvalues vanish, at least on $\{\nabla f\neq0\}$.
To sum up, the Ricci tensor of $(M^n,g)$
has at each point of $\{\nabla f\neq0\}\subset M$ the eigenvalues
$\varepsilon$ of multiplicity $2$ and $0$ of multiplicity $n-2$
respectively.
Note that both eigendistributions of the Ricci-tensor are smooth since they have
constant rank.
Furthermore, the critical set $\{\nabla f=0\}$ of $f$ must have empty
interior, otherwise the Ricci tensor would vanish identically on that interior
by (\ref{eq:nabladf=-fric}) and the fact that $0$ is not a critical value of
$f$.
But this would contradict the fact that the scalar curvature $S$ of $(M^n,g)$ is assumed to be constant and nonvanishing.
Therefore, $\Ric$ has actually $\varepsilon$ and $0$ as eigenvalues with
multiplicities $2$ and $n-2$ respectively on all of $M$.
This proves \ref{statement:valpRicci}.\\
It remains to show that, when $n=3$, both eigendistributions of
the Ricci tensor of $(M^3,g)$ are actually parallel.
Let $\eta$ be a unit eigenvector of $\Ric$ associated to the eigenvalue
$\varepsilon$ and $e_3$ be a unit eigenvector of $\Ric$ associated to the
eigenvalue $0$; since both $\Ric$-eigenvalues are constant and distinct and $\Ric$ is
smooth, $\eta$ and $e_3$ exist globally along $N_c$, no need of
analyticity.
In dimension $3$ again, because $S_{N_c}=0$ yields $\Ric_{N_c}=0$ and thus $\nabla_\nu\Ric=0$, the vector fields $\eta$ and $e_3$ can actually be defined everywhere on the regular set of $f$ using parallel transport along $\nu$-geodesics.
Moreover, because the eigenvalue $0$ of the Ricci-tensor has multiplicity $1$ on all of $M$ as we showed above, the vector field $e_3$ can be defined globally on $M$.\\
We show that
$\nabla e_3=0$, i.e. $e_3$ is parallel on the dense open subset $\{\nabla f\neq0\}$ and hence on $M$.
First, because of $\nabla_\nu\Ric=0$, $\ker(\Ric)=\R e_3$ and $|e_3|=1$, we have $\nabla_\nu e_3\in\ker(\Ric)\cap e_3^\perp=\{0\}$ i.e., $\nabla_\nu e_3=0$.
Next, following from the identity
\[0=\frac{1}{2}\nabla S=-\delta\Ric=(\nabla_\eta\Ric)\eta+(\nabla_{e_3}\Ric)e_3+\underbrace{(\nabla_\nu\Ric)}_{0}\nu,\]
we have $(\nabla_\eta\Ric)\eta=-(\nabla_{e_3}\Ric)e_3$.
Here we notice that
\[(\nabla_\eta\Ric)\eta=\varepsilon\nabla_\eta\eta-\Ric(\nabla_\eta\eta)=\varepsilon(\nabla_\eta\eta-\langle\nabla_\eta\eta,\nu\rangle\nu)\]
and, with $\nabla_{e_3}\nu=-We_3=\frac{f}{|\nabla f|}\Ric(e_3)=0$, that
\[(\nabla_{e_3}\Ric)e_3=-\Ric(\nabla_{e_3}e_3)=-\varepsilon\langle\nabla_{e_3}e_3,\eta\rangle\eta.\]
Therefore,
\be
0&=&\varepsilon\langle\nabla_\eta\eta,\eta\rangle\\
&=&\langle(\nabla_\eta\Ric)\eta,\eta\rangle\\
&=&-\langle (\nabla_{e_3}\Ric)e_3,\eta\rangle\\
&=&\varepsilon\langle\nabla_{e_3}e_3,\eta\rangle.
\ee
Since $\langle\nabla_{e_3}e_3,\nu\rangle=-\langle e_3,\nabla_{e_3}\nu\rangle=0$ and $\langle\nabla_{e_3}e_3,e_3\rangle=0$, it can be deduced that $\nabla_{e_3}e_3=0$.\\
Analogously,
\be
0&=&-\varepsilon\langle\nabla_{e_3}e_3,\eta\rangle\langle\eta,e_3\rangle\\
&=&\langle(\nabla_{e_3}\Ric)e_3,e_3\rangle\\
&=&-\langle(\nabla_\eta\Ric)\eta,e_3\rangle\\
&=&-\varepsilon\langle\nabla_\eta\eta,e_3\rangle,
\ee
so that $\langle\nabla_\eta e_3,\eta\rangle=0$.
Again, because $\langle\nabla_\eta e_3,e_3\rangle=0=\langle\nabla_\eta e_3,\nu\rangle$, it can be deduced that $\nabla_\eta e_3=0$.
To sum up, we obtain $\nabla e_3=0$ i.e., the vector field $e_3$ is parallel on $M\setminus\{\nabla f=0\}$ and hence on $M$.
As a consequence, the holonomy group of $M$ splits locally, therefore the
universal cover of $M$ is isometric to the Riemannian product
$\Sigma\times\R$ of some complete surface $\Sigma$ with $\R$.
Moreover, using formula (\ref{eq:RXnunu}) for $X=\eta$ and taking into account that $\nabla_\nu\Ric=0$, we obtain
\[R_{\eta,\nu}\nu=\left(1+\frac{Sf^2}{2|\nabla
f|^2}\right)\cdot\Ric(\eta)-\frac{f^2}{|\nabla
f|^2}\cdot\Ric^2(\eta)=\Ric(\eta)=\varepsilon\eta,\]
so that $K(\eta,\nu)=\langle
R_{\eta,\nu}\nu,\eta\rangle=\varepsilon|\eta|^2=\varepsilon$.
Therefore, the distribution $\mathrm{Span}(\eta,\nu)\to M$ integrates to a surface of constant
curvature $\varepsilon\in\{\pm1\}$.
Thus $\Sigma=S^2(\varepsilon)$, which is the
simply-connected complete surface
with curvature $\varepsilon\in\{\pm1\}$.
In case $\varepsilon=1$, the lift $\tilde{f}$ of $f$ to $\mathbb{S}^2\times\R$
is constant along the $\R$-factor and satisfies the equation
$(\nabla^{\mathbb{S}^2})^2f=-f\cdot\mathrm{Id}$, which is exactly the equation
characterizing the eigenfunctions associated to the first positive Laplace
eigenvalue \cite[Theorem A]{Obata62}.
Furthermore, the isometry group of $\mathbb{S}^2\times\R$ embeds into the product group of both isometry groups of $\mathbb{S}^2$ and $\R$ and the first factor must be trivial since $\tilde{f}$, as the restriction of a linear form from $\R^3$ onto $\mathbb{S}^2$, is not invariant under $\{\pm\mathrm{Id}\}$.
Therefore, $M$ is isometric to either $\mathbb{S}^2\times\R$ or to
$\mathbb{S}^2\times\mathbb{S}^1$ and in both cases $f$ is the trivial extension
of an eigenfunction associated to the first positive Laplace eigenvalue on
$\mathbb{S}^2$.
In case $\varepsilon=-1$, the lift $\tilde{f}$ of $f$ to $\mathbb{H}^2\times\R$
is constant along the $\R$-factor and satisfies the equation
$(\nabla^{\mathbb{H}^2})^2f=f\cdot\mathrm{Id}$, which is exactly the Tashiro
equation.
Since the isometry group of $\mathbb{H}^2\times\R$ embeds into the product group
of both isometry groups of $\mathbb{H}^2$ and $\R$ and the first factor must be
trivial since $\tilde{f}$ has no nontrivial symmetry \cite[Theorem 2
p.252]{Tashiro65}, we can deduce as above that $M$ is isometric to either
$\mathbb{H}^2\times\R$ or $\mathbb{H}^2\times\mathbb{S}^1$ and $f$ is the
trivial extension of a solution to the Tashiro equation on $\mathbb{H}^2$.
This proves statement \ref{statement:scalcstn=3} and concludes the proof of Proposition \ref{l:eqnabladf=-fricscalconstant}.
\findemo

Next we look at manifolds with \emph{harmonic} curvature tensor.
Recall that, by definition, the Riemann curvature tensor $R$ of $(M^n,g)$ is harmonic if and only if $\delta R=0$ holds on $M$.
By the first and second Bianchi identities, we have, for all $X,Y,Z\in T_xM$ at some $x\in M$:
\[(\delta R)(X,Y,Z)=(\nabla_Y\mathrm{Ric})(Z,X)-(\nabla_Z\mathrm{Ric})(Y,X).\]
As a consequence, $\delta R=0$ at some $x\in M$ is equivalent to
\[(\nabla_X\mathrm{Ric})(Y)-(\nabla_Y\mathrm{Ric})(X)=0\]
for all $X,Y\in T_xM$ i.e., to $\mathrm{Ric}$ being a \emph{Codazzi-}tensor at
$x$.
A $3$-dimensional Riemannian manifold has harmonic
curvature if and only if it is conformally flat and has constant scalar
curvature.
In dimension $n\geq4$, a Riemannian manifold has harmonic
curvature if and only if it has harmonic Weyl tensor $W$, that is, $\delta W=0$
holds on $M$, and constant scalar
curvature.
For instance, any conformally flat manifold with constant scalar
curvature has harmonic curvature tensor.
We refer to \cite[Sec. 16.4]{BesseEinstein} for more details about harmonic
curvature.

\begin{prop}\label{l:harmcurvnabladf=-fric}
Let $(M^n,g)$ be any connected Riemannian manifold carrying a nonzero smooth real-valued
function $f$ satisfying {\rm(\ref{eq:nabladf=-fric})} on $M$.
If the Riemann curvature tensor of $(M^n,g)$ is harmonic, then either $(M^n,g)$ is Ricci-flat or,
up to rescaling the metric $g$, the manifold $(M^n,g)$ is isometric to the Riemannian product $S^2(\varepsilon)\times\Sigma^{n-2}$, where
$S^2(\varepsilon)$ is the simply-connected complete surface of constant curvature
$\varepsilon\in\{\pm1\}$ and $\Sigma^{n-2}$ is a Ricci-flat manifold.
Moreover, $f$ is the trivial extension to $M$ of a solution of the Obata resp.
Tashiro equation on $\mathbb{S}^2$ (if $\varepsilon=1$) resp. $\mathbb{H}^2$ (if
$\varepsilon=-1$).
\end{prop}

{\it Proof:}
First recall that, if $\delta R=0$ holds on $M$ -- or, equivalently, if
$\mathrm{Ric}$ is a Codazzi-tensor -- then the scalar curvature $S$ of $(M^n,g)$
must be constant: given any pointwise o.n.b. $(e_j)_{1\leq j\leq n}$ of $TM$ and
$X\in TM$, we have
\be
X(S)&=&X\left(\mathrm{tr}(\mathrm{Ric})\right)\\
&=&\mathrm{tr}\left(\nabla_X\mathrm{Ric}\right)\\
&=&\sum_{j=1}^n(\nabla_X\mathrm{Ric})(e_j,e_j)\\
&=&\sum_{j=1}^n(\nabla_{e_j}\mathrm{Ric})(X,e_j)\\
&=&-(\delta\mathrm{Ric})(X)\\
&=&\frac{X(S)}{2},
\ee
so that necessarily $dS=0$ holds on $M$.
Since the scalar curvature $S$ is assumed to be non-identically vanishing, we
may assume up to rescaling $g$ that $S=2\varepsilon$ with
$\varepsilon\in\{\pm1\}$.\\
For any $s\in \mathbb{N}$, we denote by $(a_s)$ the assertion ${\rm tr}({\rm Ric}^s)=2\varepsilon^s$ and by $(b_s)$ the assertion $\delta({\rm Ric}^s)=0$.
We show that, since the Ricci-tensor is assumed to be Codazzi, both $(a_s)$ and
$(b_s)$ are true.\\
First, we have that, for every $s$, $(b_s)$ implies $(a_{s+1})$: namely, as a consequence of ${\rm Ric}(\nabla f)=\varepsilon\nabla f$ (see (\ref{eq:Ricnablaf})),
$$(\nabla_X {\rm Ric}^s)(\nabla f)=-f(\varepsilon^s{\rm Ric} X-{\rm Ric}^{s+1}X)$$
for every $X\in TM$.
This yields, in a pointwise o.n.b. $(e_j)_{1\leq j\leq n}$ of $TM$,
\begin{eqnarray*}
\delta({\rm Ric}^s)(\nabla f)&=&-\sum_{j=1}^n(\nabla_{e_j} {\rm Ric}^s)(e_j,\nabla f)\\
&=&f(\varepsilon^sS-{\rm tr}({\rm Ric}^{s+1}))\\
&=&f(2\varepsilon^{s+1}-{\rm tr}({\rm Ric}^{s+1})).
\end{eqnarray*}
Therefore, if $\delta(\mathrm{Ric}^s)=0$, then $\mathrm{tr}(\mathrm{Ric}^{s+1})=2\varepsilon^{s+1}$.
This shows the claim.
Note that here we have not used the property that $\mathrm{Ric}$ is a Codazzi-tensor.\\
Second, we have, under the condition that $\mathrm{Ric}$ is Codazzi, that $(b_{s})\Rightarrow (b_{s+1})$.
Namely assuming $(b_{s})$, assertion $(a_{s+1})$ must hold true from the previous claim.
Therefore, for every $X\in TM$,
$$\sum_{j=1}^n(\nabla_X {\rm Ric})(e_j,{\rm Ric}^{s} e_j)={\rm tr}(\nabla_X {\rm Ric}\circ {\rm Ric}^{s})=\frac{1}{s+1}X\left(({\rm tr}\left( {\rm Ric}^{s+1}\right)\right)=0.$$
Now using the fact that the Ricci-tensor is Codazzi, we compute
\begin{eqnarray*}
0&=&\sum_{j=1}^n(\nabla_X {\rm Ric})(e_j,{\rm Ric}^{s} e_j)\\
&=&\sum_{j=1}^n(\nabla_{e_j} {\rm Ric})({\rm Ric}^{s} e_j,X)\\
&=&\sum_{j=1}^n((\nabla_{e_j} {\rm Ric}^{s+1})(e_i),X)-{\rm Ric}(((\nabla_{e_j}{\rm Ric}^{s})(e_j),X)\\
&=&-(\delta\mathrm{Ric}^{s+1})(X)
\end{eqnarray*}
using again $(b_s)$.
We deduce that $(b_{s+1})$ is true. \\
Since $(a_s)$ and $(b_s)$ are satisfied for $s=1$, we deduce that they are satisfied for all $s\in\mathbb{N}$.
From the Newton identities, it can be deduced that the Ricci tensor must have pointwise the eigenvalues $\varepsilon$ and $0$, the former of multiplicity $2$ and the latter of multiplicity $n-2$.
Therefore, we get the pointwise orthogonal decomposition $TM=\ker({\rm Ric}-\varepsilon {\rm Id} )\oplus \ker({\rm Ric})$. \\

It remains to show that both eigendistributions of the Ricci-tensor are parallel.
Let $X,Y\in \ker({\rm Ric}-\varepsilon {\rm Id})$ and $Z\in \ker({\rm Ric})$.
Then the scalar product with $Y$ in the formula $(\nabla_X {\rm Ric})Z=(\nabla_Z {\rm Ric})X$ allows to get on the one hand
$$g((\nabla_X {\rm Ric})Z, Y)=-g({\rm Ric}(\nabla_X Z),Y)=-\varepsilon g(\nabla_X Z,Y),$$
and on the other hand
\begin{eqnarray*}
g((\nabla_Z {\rm Ric})X, Y)&=&\varepsilon g(\nabla_Z X,Y)-g({\rm Ric}(\nabla_Z X),Y)\\
&=&\varepsilon g(\nabla_Z X,Y)-g(\nabla_Z X,{\rm Ric}Y)\\
&=&0.
\end{eqnarray*}
Thus, we deduce that $0=g(\nabla_X Z,Y)=-g(\nabla_X Y,Z)$. Hence $\nabla_X Y\in \ker({\rm Ric}-\varepsilon\mathrm{Id})$ and therefore the distribution $\ker({\rm Ric}-\varepsilon\mathrm{Id})$ is parallel.
The same computations can be done for the distribution $\ker({\rm Ric})$.
This straightforwardly implies that both eigendistributions
$\ker(\Ric-\varepsilon\mathrm{Id})$ and $\ker(\Ric)$ are parallel and
therefore integrable and totally geodesic.
By the de Rham theorem, $M$ splits locally as the Riemannian product of a surface and
an $n-2$-dimensional sub\-ma\-ni\-fold.
Moreover, the Ricci-curvature -- which is the Gau\ss{}-curvature -- of the
surface that is pointwise tangent to the distribution
$\ker(\Ric-\varepsilon\mathrm{Id})$ is $\varepsilon$ and the sub\-ma\-ni\-fold
that is pointwise tangent to $\ker(\Ric)$ is Ricci-flat, see e.g.
\cite[Thm. 1.100]{BesseEinstein}.
Therefore the universal cover of $M$ is isometric to the Riemannian product
$S^2(\varepsilon)\times\tilde{\Sigma}$ of the simply-connected complete surface
with curvature $\varepsilon\in\{-1,0,1\}$ with some simply-connected Ricci-flat
manifold $\tilde{\Sigma}$.
The rest of the proof is analogous to that of Proposition \ref{l:eqnabladf=-fricscalconstant}.\ref{statement:scalcstn=3}.
This concludes the proof of
Proposition \ref{l:harmcurvnabladf=-fric}.
\findemo

\section{Examples in warped product form}\label{s:exwp}

We look for examples of warped products $(M,g):=(M_1\times
M_2,g_1\oplus\varphi^2 g_2)$ for some smooth positive function $\varphi$ on 
$M_1$, where $(M_1,g_1)$ and $(M_2,g_2)$ are connected Riemannian manifolds.
We make the ansatz $f(x_1,x_2):=f_1(x_1)f_2(x_2)$ for all $(x_1,x_2)\in M$ where $f_1$ and $f_2$ are smooth real-valued functions on $M_1$ and $M_2$ respectively.
We 
look for necessary and sufficient conditions for $f$ to satisfy (\ref{eq:nabladf=-fric}) on 
$(M,g)$.

\begin{prop}\label{p:eqhessianricciwp01}
Let $(M^n,g):=(M_1^{n_1}\times M_2^{n_2},g_1\oplus\varphi^2g_2)$ be a connected Riemannian warped product, where $\varphi\in C^\infty(M_1,\R_+^\times)$.
For any two functions $f_i\in C^\infty(M_i,\R)$, $i=1,2$, let $f:=\pi_1^*f_1\cdot\pi_2^*f_2$ i.e., $f(x_1,x_2)=f_1(x_1)f_2(x_2)$ for all $(x_1,x_2)\in M$.
Then $f$ solves $\nabla^2 f=-f\cdot\Ric$ on $(M,g)$ if and only if one of the following occurs:
\begin{enumerate}[(a)]
\item The function $\frac{f_1}{\varphi}$ is constant on $M_1$, in which case it can be assumed up to re\-sca\-ling $f$ that $f_1=\varphi$.
Then $\mu_1(f_1):=(n_2-2)|\nabla^{M_1}f_1|_1^2-f_1\Delta^{M_1}f_1$
is constant on $M_1$ and $f_1,f_2$ solve
\begin{eqnarray}\label{eq:nabla2f1casef1=phi}(n_2-1)(\nabla^{M_1})^2f_1&=&f_1\cdot\Ric_{M_1}\\\label{eq:nabla2f2casef1=phi}(\nabla^{M_2})^2f_2&=&f_2\cdot\left(\mu_1(f_1)\mathrm{Id}_{TM_2}-\Ric_{M_2}\right)
\end{eqnarray}
respectively.
\item The function $f_2$ is constant on $M_2$, in which case $f_1$ solves
\begin{equation}\label{eq:nabla2f1casef2constant}
(\nabla^{M_1})^2f_1=-f_1\cdot\left(\Ric_{M_1}-
\frac{n_2}{\varphi}(\nabla^{M_1})^2\varphi\right)
\end{equation}
on $M_1$, the function $-\frac{\varphi}{f_1}g_1(\nabla^{M_1}f_1,\nabla^{M_1}
\varphi)+(n_2-1)|\nabla^{M_1}\varphi|_1^2-\varphi\Delta^{M_1}\varphi$ is
constant on each connected component of $M_1\setminus f_1^{-1}(\{0\})$ and the
manifold $(M_2,g_2)$ is Einstein with scalar curvature equal to
\[n_2\left(-\frac{\varphi}{f_1}g_1(\nabla^{M_1}f_1,\nabla^{M_1}
\varphi)+(n_2-1)|\nabla^{M_1}\varphi|_1^2-\varphi\Delta^{M_1}\varphi\right).\]
\end{enumerate}
\end{prop}

{\it Proof:} First, we have $\nabla f=f_2\nabla f_1+f_1\nabla
f_2=f_2\nabla^{M_1}f_1+\frac{f_1}{\varphi^2}\nabla^{M_2}f_2$, where 
$\nabla^{M_i}f_i$ denotes the $g_i$-gradient of $f_i$ on $(M_i,g_i)$.
Recall Koszul's formula, valid for any tangent vector fields $X,Y,Z$ on some Riemannian manifold $(M,g)$:
\begin{eqnarray}\label{eq:Koszulformula} 
\nonumber g(\nabla_X
Y,Z)&=&\frac{1}{2}\Big\{X(g(Y,Z))+Y(g(Z,X))-Z(g(X,Y))\\
&&\phantom{\frac{1}{2}\Big\{}+g([X,Y],Z)-g([Y,Z],X)+g([Z,X],Y)\Big\}.
\end{eqnarray}
It can be deduced from (\ref{eq:Koszulformula}) that, for any 
$X_i,Y_i,Z_i\in\Gamma(\pi_i^* TM_i)$, we have 
\begin{eqnarray}
\label{eq:nabla11}\nabla_{X_1}Y_1&=&\nabla_{X_1}^{M_1}Y_1\\
\label{eq:nabla12}\nabla_{X_1}Y_2&=&\partial_{X_1}Y_2+\frac{X_1(\varphi)}{
\varphi}Y_2\\
\label{eq:nabla21}\nabla_{X_2}Y_1&=&\partial_{X_2}Y_1+\frac{Y_1(\varphi)}{
\varphi}X_2\\
\label{eq:nabla22}\nabla_{X_2}Y_2&=&\nabla_{X_2}^{M_2}Y_2-\frac{1}{\varphi}g(X_2
,Y_2)\nabla^{M_1} \varphi.
\end{eqnarray}
As a first consequence,
\begin{eqnarray}\label{eq:HessianfX1}
\nonumber\nabla_{X_1}^2f&=&f_2\nabla_{X_1}\nabla^{M_1}f_1+\frac{X_1(f_1)\varphi^2-2f_1X_1
(\varphi)\varphi} 
{\varphi^4}\nabla^{M_2}f_2+\frac{f_1}{\varphi^2}\nabla_{X_1}\nabla^{M_2}
f_2\\
\nonumber&=&f_2(\nabla^{M_1})_{X_1}^2f_1+\frac{X_1(f_1)\varphi-2f_1X_1
(\varphi)} 
{\varphi^3}\nabla^{M_2}f_2\\
\nonumber&&+\frac{f_1}{\varphi^2}\left(\underbrace{\partial_{X_1}\nabla^{M_2}f_2}_{
0}+\frac { X_1(\varphi) } { \varphi}\nabla^{M_2}f_2\right)\\
\nonumber&=&f_2(\nabla^{M_1})_{X_1}^2f_1+\frac{X_1(f_1)\varphi-f_1X_1
(\varphi)} 
{\varphi^3}\nabla^{M_2}f_2\\
&=&f_2(\nabla^{M_1})_{X_1}^2f_1+\frac{1}{\varphi}X_1(\frac{f_1}{
\varphi} )\nabla^{M_2} f_2.
\end{eqnarray}
Similarly,
\begin{eqnarray}\label{eq:HessianfX2}
\nonumber\nabla_{X_2}^2f&=&X_2(f_2)\nabla^{M_1}f_1+f_2\nabla_{X_2}\nabla^{M_1}f_1+\frac
{f_1}{\varphi^2}\nabla_{X_2}\nabla^{M_2}f_2\\
\nonumber&=&X_2(f_2)\nabla^{M_1}f_1+f_2\left(\underbrace{\partial_{X_2}\nabla^{M_1}
f_1}_{0}+\frac { g_1(\nabla^{M_1}f_1,\nabla^{M_1}\varphi)}{\varphi}X_2\right)\\
\nonumber&&+\frac{f_1}{\varphi^2}\left((\nabla^{M_2})_{X_2}^2f_2-\frac{1}{\varphi}g(X_2,\nabla^{
M_2}f_2)\nabla^{M_1}\varphi\right)\\
\nonumber&=&\frac{f_1}{\varphi^2}(\nabla^{M_2})_{X_2}^2f_2+X_2(f_2)\left(\nabla^{M_1}f_1-\frac{f_1}{\varphi}\nabla^{M_1}\varphi\right)\\
\nonumber&&+\frac{f_2}{\varphi}g_1(\nabla^{M_1}f_1,\nabla^{M_1}\varphi)X_2\\
&=&\frac{f_1}{\varphi^2}(\nabla^{M_2})_{X_2}^2f_2+X_2(f_2)\varphi\nabla^{M_1}\left(\frac{f_1}{\varphi}\right)+\frac{f_2}{\varphi}g_1(\nabla^{M_1}f_1,\nabla^{M_1}\varphi)X_2.
\end{eqnarray}
Independently, by \cite[Prop. 9.106]{BesseEinstein}, we have 
\begin{eqnarray}
\label{eq:RicM1wp}\Ric(X_1)&=&\Ric_{M_1}(X_1)-\frac{n_2}{\varphi}(\nabla^{M_1})_{X_1}^2\varphi\\
\label{eq:RicM2wp}\Ric(X_2)&=&\frac{1}{\varphi^2}\Ric_{M_2}(X_2)+\left(\frac{\Delta^{M_1}\varphi
}{\varphi}-(n_2-1)\frac{|\nabla^{M_1}\varphi|_1^2}{\varphi^2}\right)X_2
\end{eqnarray}
Therefore, $f$ satisfies (\ref{eq:nabladf=-fric}) on $(M,g)$ if and only if the following system of equations holds, for all $(X_1,X_2)\in TM$:
\[\left\{\begin{array}{ll}l_1(X_1)&=r_1(X_1)\\ l_2(X_2)&=r_2(X_2)\end{array}\right.,\]
where 
\begin{eqnarray*}
l_1(X_1)&=&f_2(\nabla^{M_1})_{X_1}^2f_1+\frac{1}{
\varphi}X_1(\frac{f_1}{
\varphi} )\nabla^{M_2} f_2\\
r_1(X_1)&=&-f_1f_2\cdot\left(\Ric_{M_1}(X_1)-\frac{n_2}{\varphi}(\nabla^{M_1})_{X_1}^2\varphi\right)\\
l_2(X_2)&=&\frac{f_1}{\varphi^2}(\nabla^{M_2})_{X_2}^2f_2+X_2(f_2)\varphi\nabla^{M_1}\left(\frac{f_1}{\varphi}\right)+\frac{f_2}{\varphi}g_1(\nabla^{M_1}f_1,\nabla^{M_1}\varphi)X_2\\
r_2(X_2)&=&-f_1f_2\cdot\left(\frac{1}{\varphi^2}\Ric_{M_2}(X_2)+\left(\frac{\Delta^{M_1}\varphi
}{\varphi}-(n_2-1)\frac{|\nabla^{M_1}\varphi|_1^2}{\varphi^2}\right)X_2\right).
\end{eqnarray*}
\small
\normalsize
Both equations imply that $d\left(\frac{f_1}{\varphi}\right)\otimes df_2=0$, that is, that $\frac{f_1}{\varphi}$ is constant on $M_1$ or $f_2$ is constant on $M_2$.\\
{\bf Case $\frac{f_1}{\varphi}$ is constant on $M_1$:} We may assume, up to re\-sca\-ling $f_2$ and hence $f$, that $f_1=\varphi$ holds on $M_1$.
The above system of equations becomes equivalent to the following:
\[\left\{\begin{array}{ll}(\nabla^{M_1})_{X_1}^2f_1&=-f_1\cdot\left(\Ric_{M_1}(X_1)-\frac{n_2}{f_1}(\nabla^{M_1})_{X_1}^2f_1\right)\\\frac{1}{f_1}(\nabla^{M_2})_{X_2}^2f_2+\frac{f_2|\nabla^{M_1}f_1|_1^2}{f_1}X_2&=-f_1f_2\cdot T_2(X_2)\end{array}\right.,\]
where $T_2(X_2):=\left(\frac{1}{f_1^2}\Ric_{M_2}(X_2)+\left(\frac{\Delta^{M_1}f_1
}{f_1}-(n_2-1)\frac{|\nabla^{M_1}f_1|_1^2}{f_1^2}\right)X_2\right)$.
Thus
\[\left\{\begin{array}{ll}(1-n_2)(\nabla^{M_1})^2f_1&=-f_1\cdot\Ric_{M_1}\\(\nabla^{M_2})^2f_2&=-f_2\cdot\Ric_{M_2}+f_2\left((n_2-2)|\nabla^{M_1}f_1|_1^2-f_1\Delta^{M_1}f_1\right)\mathrm{Id}_{TM_2}\end{array}\right..\]
Since $f_2$ is assumed to be non-identically vanishing and the second identity above only depends on $M_2$, the factor $\mu_1(f_1):=(n_2-2)|\nabla^{M_1}f_1|_1^2-f_1\Delta^{M_1}f_1$ must be constant on $M_1$.
Actually we shall see later that, when $n_2\geq2$, this already follows from the equation for $f_1$.\\
Therefore, in case $f_1=\varphi$, equation (\ref{eq:nabladf=-fric}) for $f:=\pi_1^*f_1\cdot\pi_2^*f_2$ is equivalent to the function $(n_2-2)|\nabla^{M_1}f_1|_1^2-f_1\Delta^{M_1}f_1=\mu_1(f_1)$ begin constant on $M_1$ and
\[
\left\{\begin{array}{ll}(n_2-1)(\nabla^{M_1})^2f_1&=f_1\cdot\Ric_{M_1}\\(\nabla^{M_2})^2f_2&=f_2\cdot\left(\mu_1(f_1)\mathrm{Id}_{TM_2}-\Ric_{M_2}\right)\end{array}\right..
\]
{\bf Case $f_2$ is constant on $M_2$:} 
Then $\nabla^2 f=-f\cdot\Ric$ on $(M,g)$ is equivalent to the system
\[\left\{\begin{array}{ll}(\nabla^{M_1})^2f_1&=-f_1\cdot\left(\Ric_{M_1}-\frac{
n_2}{\varphi}(\nabla^{M_1})^2\varphi\right)\\ 
\frac{g_1(\nabla^{M_1}f_1,\nabla^{M_1}\varphi)}{\varphi}\mathrm{Id}_{TM_2}
&=-f_1\cdot\left(\frac{1}{\varphi^2}\Ric_{M_2}+(\frac{\Delta^{M_1}\varphi}{
\varphi}-(n_2-1)\frac{|\nabla^{M_1}\varphi|_1^2}{\varphi^2})\mathrm{Id}_{TM_2}
\right)\end{array}\right.,\]
that is, assuming $f_1$ not to vanish identically on $M_1$,
\[\left\{\begin{array}{ll}(\nabla^{M_1})^2f_1&=-f_1\cdot\left(\Ric_{M_1}-
\frac{n_2}{\varphi}(\nabla^{M_1})^2\varphi\right)\\ 
\Ric_{M_2}
&=\left(-\frac{\varphi}{f_1}g_1(\nabla^{M_1}f_1,\nabla^{M_1}
\varphi)+(n_2-1)|\nabla^{M_1}\varphi|_1^2-\varphi\Delta^{M_1}\varphi
\right)\cdot\mathrm{Id}_{TM_2}\end{array}\right.,\]
the se\-cond equation holding on the dense open subset $M_1\setminus f^{-1}(\{0\})$.
The se\-cond of both above identities implies that the quantity
\[\mu_1':=\left(-\frac{\varphi}{f_1}g_1(\nabla^{M_1}f_1,\nabla^{M_1}
\varphi)+(n_2-1)|\nabla^{M_1}\varphi|_1^2-\varphi\Delta^{M_1}\varphi\right)\]
is constant on $M_1$ and that $M_2$ is Einstein with constant scalar curvature equal to $n_2\mu_1'$, whatever $n_2$ is.
This concludes the proof of Proposition \ref{p:eqhessianricciwp01}.
\findemo

Now we look at (\ref{eq:nabladf=-fric}) on Riemannian products, where $f$ is not assumed to be in product form:

\begin{prop}\label{l:wpcasephi=1}
Let $(M^n,g)=(M_1\times M_2,g_1\oplus g_2)$ for some connected Riemannian manifolds $(M_1,g_1)$ and $(M_2,g_2)$.
Assume $M$ to be non Ricci-flat i.e., that $\Ric_{M_1}\neq0$ or $\Ric_{M_2}\neq0$.
W.l.o.g. let $\Ric_{M_2}\neq0$.
Then a function $f\in C^\infty(M,\R)\setminus\{0\}$ satisfies {\rm (\ref{eq:nabladf=-fric})} on $(M^n,g)$ if and only if $\Ric_{M_1}=0$, the function $f$ only depends on $M_2$ and satisfies  {\rm (\ref{eq:nabladf=-fric})} on $(M_2,g_2)$.
As a consequence, the map $W(M_2,g_2)\lra W(M,g)$ extending a function trivially on the $M_1$-factor, is an isomorphism.
\end{prop}

{\it Proof:}
First, we split pointwise $\nabla f=\nabla^{M_1}f+\nabla^{M_2}f$ according to the $g$-orthogonal splitting $T_{(x_1,x_2)}M=T_{x_1}M_1\oplus T_{x_2}M_2$, for all $(x_1,x_2)\in M$.
Using formulae (\ref{eq:nabla11}) -- (\ref{eq:nabla22}) and $\varphi=1$, it can be deduced that, for all $X_1\in TM_1$,
\begin{eqnarray*}
\nonumber\nabla_{X_1}^2f&=&\nabla_{X_1}(\nabla^{M_1}f)+\nabla_{X_1}(\nabla^{M_2}f)\\
&=&\nabla_{X_1}^{M_1}(\nabla^{M_1}f)+\partial_{X_1}(\nabla^{M_2}f)
\end{eqnarray*}
and similarly, for all $X_2\in TM_2$,
\begin{eqnarray*}
\nonumber\nabla_{X_2}^2f&=&\nabla_{X_2}(\nabla^{M_1}f)+\nabla_{X_2}(\nabla^{M_2}f)\\
&=&\partial_{X_2}(\nabla^{M_1}f)+\nabla_{X_2}^{M_2}(\nabla^{M_2}f).
\end{eqnarray*}
By (\ref{eq:RicM1wp}) and (\ref{eq:RicM2wp}), we obtain that $f$ satisfies (\ref{eq:nabladf=-fric}) on $(M^n,g)$ if and only if, for all $(X_1,X_2)\in TM_1\oplus TM_2$,
\begin{eqnarray}
\nabla_{X_1}^{M_1}(\nabla^{M_1}f)+\partial_{X_1}(\nabla^{M_2}f)\label{eq:nablaX1nablaM1f}&=&-f\Ric_{M_1}(X_1)\\
\partial_{X_2}(\nabla^{M_1}f)+\nabla_{X_2}^{M_2}(\nabla^{M_2}f)\label{eq:nablaX2nablaM1f}&=&-f\Ric_{M_2}(X_2).
\end{eqnarray}
It can be deduced that both $\partial_{X_1}(\nabla^{M_2}f)=0$ and $\partial_{X_2}(\nabla^{M_1}f)=0$, for all $(X_1,X_2)\in TM_1\oplus TM_2$.
But the first identity is equivalent to the existence of functions $a_1\in C^\infty(M_1,\R)$ and $a_2\in C^\infty(M_2,\R)$ such that $f(x_1,x_2)=a_1(x_1)+a_2(x_2)$ for all $(x_1,x_2)\in M$.
Then the second identity is trivial and (\ref{eq:nablaX2nablaM1f}) is equivalent to
\[(\nabla^{M_2})^2a_2{}_{|_{x_2}}=-(a_1(x_1)+a_2(x_2))\Ric_{M_2}{}_{|_{x_2}}\]
for all $(x_1,x_2)\in M_1\times M_2$.
But since the l.h.s. of the preceding inequality does not depend on $M_1$ and because of $\Ric_{M_2}\neq0$, this implies $a_1$ is constant on $M_1$, therefore $a_1+a_2\in C^\infty(M_2,\R)$ satisfies (\ref{eq:nabladf=-fric}) on $(M_2,g_2)$.
But then (\ref{eq:nablaX1nablaM1f}) together with the assumption $f\neq0$ forces $\Ric_{M_1}=0$: choose a point $x_2\in M_2$ where $f(x_2)\neq0$.
This concludes the proof.
\findemo

\color{black}

Next we look for examples and partial classifications results for identities 
(\ref{eq:nabla2f1casef1=phi}) and (\ref{eq:nabla2f2casef1=phi}), which 
correspond to the case $f_1=\varphi$.
 An obvious case is when $f_1=\varphi$ are constant (and
nonvanishing) on $M_1$.
Then $(M_1,g_1)$ must be Ricci-flat, $f(x_1,x_2)=f_2(x_2)$ for all 
$(x_1,x_2)\in M$ and, because of $\mu_1(f_1)=0$ then, the function $f_2$ must 
satisfy (\ref{eq:nabladf=-fric}) on $(M_2,g_2)$.
This is actually a consequence of Proposition \ref{l:wpcasephi=1} above.
Therefore we obtain an already known example in that case, see 
introduction. \color{black}

\begin{prop}\label{p:hessianricciwarpedproducts}
Let $(M^n,g)$ be any connected Riemannian manifold.
\begin{enumerate}
\item\label{statement:f1} Assume there exists an $f\in C^\infty(M,\R_+^\times)$ solving
$\nabla^2f=\frac{f}{m-1}\cdot\Ric$ on $M$ for some integer $m\geq2$.
Then $\mu_1(f):=(m-2)|\nabla f|^2+\frac{f^2S}{m-1}=(m-2)|\nabla f|^2-f\Delta f$
is constant on $M$ and, if $m>2$, then
$\Ric(\nabla f)=-\frac{1}{(m-2)(m-1)}\nabla(f^2S)$, where $S$ is the scalar
curvature of $(M,g)$.
Moreover, if $m>2$, then $f$ defines a $(0,n+m-1)$-Einstein metric on $(M,g)$.
\item\label{statement:f2} Assume there exists an $f\in C^\infty(M,\R_+^\times)$ solving
$\nabla^2f=f\cdot(\mu\mathrm{Id}-\Ric)$ on $M$ for some $\mu\in\R$.
Then
$\Ric(\nabla
f)=-\frac{(n-1)\mu}{2}\nabla f+\frac
{f}{4}\nabla S+\frac{S}{2}\nabla f$ and
$\mu_2(f):=2|\nabla f|^2+f^2(S-(n+1)\mu)=2|\nabla
f|^2+f\Delta f-\mu f^2$ is constant on $M$.
\item\label{statement:exampleswp} In case $(M^n,g)=(M_1^{n_1}\times M_2^{n_2},g_1\oplus f_1^2g_2)$ for some
$f_1\in C^\infty(M_1,\R_+^\times)$ and $f:=\pi_1^*f_1\cdot\pi_2^*f_2$ for some
$f_2\in C^\infty(M_2,\R)$, there are, for each $n_1,n_2\geq1$ examples of
$(M_i,g_i,f_i)$ for which $f$ solves {\rm (\ref{eq:nabladf=-fric})}.
\item\label{statement:Mclosed} If $(M^n,g)$ is closed and $f\in
C^\infty(M,\R_+^\times)$ is such that $\mu_1(f):=k|\nabla f|^2-f\Delta f$ is
constant for some $k\in\R$, then $f$ must be
constant on $M$ and therefore $\mu_1(f)$ must vanish.
As a consequence, if there exists a nonzero $f\in C^\infty(M,\R_+^\times)$ solving
$\nabla^2f=\frac{f}{m-1}\cdot\Ric$ on some closed $M$ and for some integer $m\geq2$, then $f$ must be constant and therefore  $M$ must be Ricci-flat.
\end{enumerate}
\end{prop}

{\it Proof:}
We first look at equation
\begin{equation}\label{eq:nabla2f1casef1=phin2>1}
\nabla^2f=\frac{f}{m-1}\cdot\Ric
\end{equation}
on $M$, for some integer $m\geq2$.
We first derive a few identities following from
(\ref{eq:nabla2f1casef1=phin2>1}), see e.g. \cite[Lemma 4]{KimKim03}.
We write down the proof for the sake of completeness.
Namely, by (\ref{eq:Weitzenboeck}), we know that
\begin{eqnarray*}
\delta\left(\nabla^2f\right)&=&\Delta(\nabla f)-\Ric(\nabla f)\\
&=&\nabla(\Delta f)-\Ric(\nabla f)\\
&=&-\frac{1}{m-1}\nabla(fS)-\Ric(\nabla f)\\
&=&-\frac{1}{m-1}\left(S\nabla f+f\nabla S\right)-\Ric(\nabla f),
\end{eqnarray*}
where, as above, $S:=\mathrm{tr}(\Ric)$ is the scalar curvature of $(M,g)$ and where we have used $\Delta f=-\frac{fS}{m-1}$ tracing (\ref{eq:nabla2f1casef1=phin2>1}).
But (\ref{eq:nabla2f1casef1=phin2>1}) also yields
\begin{eqnarray*}
\delta\left(\nabla^2f\right)&=&\frac{1}{m-1}\left(-\Ric(\nabla f)+f\delta(\Ric)\right)\\
&=&\frac{1}{m-1}\left(-\Ric(\nabla f)-\frac{f}{2}\nabla S\right),
\end{eqnarray*}
so that, bringing both identities for $\delta\left(\nabla^2f\right)$ together, we deduce that
\begin{eqnarray}
\nonumber\frac{m-2}{m-1}\cdot\Ric(\nabla f)&=&-\frac{1}{m-1}\left(S\nabla f+\frac{f}{2}\nabla S\right)\\
\label{eq:Ricnablaf1n2>1}&=&-\frac{1}{2(m-1)f}\cdot\nabla(f^2S).
\end{eqnarray}
In case $m=2$, we deduce that $\nabla(f^2S)=0$ i.e., that $f^2S=-f\Delta f$ is constant on $M$.
In case $m>2$, we deduce that
\begin{equation}\label{eq:Ricnablaf1n2>2}
\Ric(\nabla f)=-\frac{1}{2(m-2)f}\cdot\nabla(f^2S).
\end{equation}
Still when $m>2$, it follows that
\begin{eqnarray*}
\nabla\left(|\nabla f|^2\right)&=&2(\nabla)_{\nabla f}^2f\\
&=&\frac{2f}{m-1}\Ric(\nabla f)\\
&\bui{=}{\rm(\ref{eq:Ricnablaf1n2>2})}&-\frac{1}{(m-2)(m-1)}\cdot\nabla(f^2S).
\end{eqnarray*}
Therefore, $\mu_1(f):=(m-2)|\nabla f|^2+\frac{1}{m-1}f^2S=(m-2)|\nabla f|^2-f\Delta f$ is constant on $M$.
Note that this is also the case when $m=2$ by the above remark.
Note also that, when $m>2$, identity (\ref{eq:nabla2f1casef1=phin2>1}) defines
a so-called $(0,n+m-1)$-Einstein metric on $(M,g)$ according to
\cite{HePetersenWylie11102455,HePetersenWylie11102456}.
By \cite[Theorem 2.2]{CaseShuWei2011}, the existence of such a positive $f$
is equivalent to the warped product $(M\times F,g\oplus f^2 g_F)$ being
Ricci-flat, where $(F,g_F)$ is an Einstein manifold of dimension $m-1$ and
with $\Ric_F=\mu_1\cdot\mathrm{Id}$, the constant $\mu_1$ being given
by
$\mu_1=(m-2)|\nabla f|^2-f\Delta f=(m-2)|\nabla f|^2+\frac{f^2S}{m-1}$, which is exactly the constant $\mu_1(f)$ described above, see also
\cite[Cor. 3]{KimKim03}.
This statement remains true when $m=2$ and
$\Delta f=0$ (or equivalently $\mu_1(f)=0$).
This shows statement \ref{statement:f1}.\\
Next we look at
\begin{equation}\label{eq:f2n2geq2}
\nabla^2f=f\cdot(\mu\mathrm{Id}-\Ric)
\end{equation}
on $M^n$ for some $\mu\in\R$ and $n\geq2$.
First and as before, a few identities can be deduced from
(\ref{eq:f2n2geq2}).
Namely, by (\ref{eq:Weitzenboeck}), we know that
\begin{eqnarray*}
\delta\left(\nabla^2f\right)&=&\Delta(\nabla
f)-\Ric(\nabla f)\\
&=&\nabla(\Delta f)-\Ric(\nabla f)\\
&=&\nabla(f(S-n\mu))-\Ric(\nabla f)\\
&=&(S-n\mu)\nabla f+f\nabla S-\Ric(\nabla f),
\end{eqnarray*}
where we have used
$\Delta f=f(S-n\mu)$ tracing (\ref{eq:f2n2geq2}).
But (\ref{eq:f2n2geq2}) also yields
\begin{eqnarray*}
\delta\left(\nabla^2f\right)&=&-(\mu\nabla f-\Ric
(\nabla f))+f\delta\left(\mu\mathrm{Id}-\Ric\right)\\
&=&-\mu\nabla f+\Ric
(\nabla f)+\frac{f}{2}\nabla S
\end{eqnarray*}
so that, bringing both identities for
$\delta\left(\nabla^2f\right)$ together, we deduce that
\begin{eqnarray}
\label{eq:Ricnablaf2n2>1}\Ric(\nabla f)&=&-\frac{n-1}{2}
\mu\nabla f+\frac
{f}{4}\nabla S+\frac{S}{2}\nabla f.
\end{eqnarray}
It follows that
\begin{eqnarray*}
\nabla\left(|\nabla f|^2\right)&=&2\nabla_{\nabla f}
^2f\\
&=&2f(\mu\nabla f-\Ric(\nabla f))\\
&\bui{=}{\rm(\ref{eq:Ricnablaf2n2>1})}&f\left((n+1)\mu\nabla f-
\frac
{f}{2}\nabla S-S\nabla f\right)\\
&=&\frac{n+1}{2}\mu\nabla(f^2)-\frac{1}{2}\nabla(f^2S)\\
&=&\frac{1}{2}\nabla\left((n+1)\mu f^2-f^2S\right).
\end{eqnarray*}
Therefore,
$\mu_2(f):=2|\nabla f|^2+f^2(S-(n+1)\mu)=2|\nabla f|^2+f\Delta f-\mu f^2$ is constant on $M$.
This proves statement \ref{statement:f2}.\\
As for statement \ref{statement:exampleswp}, we look at different cases according to the values of $n_2$ and $n_1$.\\
{\bf Case $n_2=1$:} \color{black} Then (\ref{eq:nabla2f1casef1=phi}) is equivalent to $M_1$ being Ricci-flat.
Together with $f_1\Delta^{M_1}f_1+|\nabla^{M_1}f_1|_1^2=-\mu_1(f_1)$ being constant by Proposition \ref{p:eqhessianricciwp01}, identity (\ref{eq:nabla2f2casef1=phi}) is equivalent to $f_2''=\mu_1(f_1)f_2$.
Whatever the sign of $\mu_1(f_1)$, there exists a solution $f_2$ to that se\-cond-order linear ODE on $\R$, which is periodic (and hence can be pulled down on a circle of suitable radius) if and only if $\mu_1(f_1)<0$.
As for $f_1$, a trivial family of examples in each dimension $n_1$ may be produced as follows.
When $n_1=1$, the function $f_1$ solves the ODE 
$-f_1f_1''+(f_1')^2=-\mu_1(f_1)$, whose general solution is
\[f_1(t)=\left\{\begin{array}{ll}a_1(t)&
\textrm{if }\mu_1(f_1)>0\\
b_1(t)&\textrm{if }\mu_1(f_1)=0\\ 
c_1(t),d_1(t), 
e_1(t)&\textrm{if 
}\mu_1(f_1)<0\end{array}\right.\]
where 
\begin{eqnarray*}
a_1(t)&:=&A\cosh(A^{-1}\sqrt{\mu_1(f_1)}
t+\phi)\\
b_1(t)&:=&A e^{\phi t}\\
c_1(t)&:=&A\cos(A^{-1}\sqrt{-\mu_1(f_1)}t+\phi)\\
d_1(t)&:=&\pm\sqrt{-\mu_1(f_1)}t+\phi\\
e_1(t)&:=&A\sinh(A^{-1}\sqrt{-\mu_1(f_1)}t+\phi)
\end{eqnarray*}
for real arbitrary constants $A,\phi$ with w.l.o.g. $A>0$ (remember that 
$f_1=\varphi>0$ by assumption). 
\color{black}
Note that all solutions are defined on $\R$ but that, in case $\mu_1(f_1)<0$, the function $f_1$ must change sign somewhere, which makes the solution $f_1$ only local then.
Moreover, in case $\mu_1(f_1)\geq0$, the solution $f_1$ -- though positive on $\R$ -- is not periodic and therefore cannot be pulled down on an $\mathbb{S}^1$.
Obviously, each of the above $f_1$'s can be trivially extended constantly in the 
other variables on $\R^{n_1}$ for every $n_1\geq1$.

It is important to note that, in the cases where $f_1>0$ on $\R$, 
corresponding to $\mu_1(f_1)\geq0$ as we have seen above, the induced metric 
$ds^2\oplus f_1(s)^2dt^2$ on $\R^2$ is the hyperbolic one, for which we can 
anyway describe $W(M,g)$ explicitly.
\color{black}

{\bf Case $n_2>1$:}
When $n_2\geq2$ and $n_1=1$, equation (\ref{eq:nabla2f1casef1=phi}) reduces to $f_1''=0$ on $M_1$, which has no positive solution on $M_1$ unless $f_1$ is constant or $M_1$ is a strict open subinterval of $\R$.\\
When $n_2\geq2$ and $n_1=2$, equation (\ref{eq:nabla2f1casef1=phi}) is equivalent to $(\nabla^{M_1})^2f_1=f_1\phi_1\cdot\mathrm{Id}_{TM_1}$, where $\phi_1:=\frac{S_1}{2(n_2-1)}$.
But by \cite[Sec. 2]{Tashiro65}, this implies that, on any open subset 
where $f_1$ has no critical point, $(M_1^2,g_1)$ is locally isometric to 
$(\R^2,dt^2\oplus\rho(t)^2ds^2)$, where $\rho:=\frac{u'}{u'(0)}$ and $u$ is 
$f_1$ along the flow of its normalised gradient $\nu:=\frac{\nabla^{M_1}
f_1}{|\nabla^{M_1} f_1|_1}$.
Moreover, along any integral curve $\gamma$ of $\nu$, which is a geodesic of 
$(M_1,g_1)$ because of $\nabla^{M_1} f_1$ being a pointwise eigenvector of 
$(\nabla^{M_1})^2 f_1$, the function $u$ must satisfy the following 
se\-cond-order ODE: for any $t$ in some nonempty open interval,
\begin{eqnarray*}
u''(t)&=&g_1((\nabla^{M_1})_{\dot{\gamma}(t)}^2f_1,\dot{\gamma}
(t))\\
&=&\frac{(f_1S_1)\circ\gamma(t)}{2(n_2-1)}\\               
&=&\left(\frac{\mu_1(f_1)}{2f_1}-\frac{n_2-2}{2}\cdot\frac{|\nabla^{
M_1}f_1|_1^2}{f_1}\right)\circ\gamma(t)\\
&=&-\frac{n_2-2}{2u(t)}u'(t)^2+\frac{\mu_1(f_1)}{2u(t)},
\end{eqnarray*}
that is, 
\begin{equation}\label{eq:ODEf1n1=2}
u''\cdot u+\frac{n_2-2}{2}(u')^2=\frac{\mu_1(f_1)}{2}.
\end{equation}
In the first special case where $\mu_1(f_1)=0$, the general form of the 
solution $u$ to (\ref{eq:ODEf1n1=2}) is $u(t)=(at+b)^{\frac{2}{n_2}}$ for real 
constants $a,b$ with $a\neq0$; assuming $a$ and $b$ to be positive, the maximal 
existence interval for $u$ is $[-\frac{b}{a},\infty)$, in particular no complete 
$M_1$ can exist unless $f_1$ has critical points.\\
In the second special case where $n_2=2$, the second-order ODE 
(\ref{eq:ODEf1n1=2}) may be reduced to the first-order one 
\[u'=\sqrt{\mu_1(f_1)\ln(u)+C}\]
for some real constant $C$.
Note that this implies that $u$ is constant when $n_2=2$ and $\mu_1(f_1)=0$.
If $\mu_1(f_1)>0$, the maximal existence interval for $u$ is of the form 
$]a,\infty[$, whereas if $\mu_1(f_1)<0$, that interval is of the form 
$]-\infty,a[$ for some real $a$.\\
Conversely, let us assume $u$ to be any positive solution with w.l.o.g. 
positive first derivative of (\ref{eq:ODEf1n1=2}) 
on some open interval $I$ about $0$.
Consider the warped product 
$(M_1,g_1):=(I\times\Sigma,dt^2\oplus\varphi(t)^2ds^2)$ for $\Sigma=\R$ or $\mathbb{S}^1$, where $\varphi(t):=\frac{u'(t)}{u'(0)}$.
Let $f(t,s):=u(t)$ for all $(t,s)\in M_1$.
The above formulae (\ref{eq:HessianfX1}) and (\ref{eq:HessianfX2}) for the Hessian of $f$ simplify to $\nabla_{\partial_t}^2f=u''\cdot\partial_t$ and $\nabla_{\partial_s}^2f=\frac{u'\varphi'}{\varphi}\cdot\partial_s$.
The identities (\ref{eq:RicM1wp}) and (\ref{eq:RicM2wp}) become $\Ric=-\frac{\varphi''}{\varphi}\cdot\mathrm{Id}_{TM}$.
Taking into account that $\varphi=\frac{u'}{u'(0)}$, we have $\frac{u'\varphi'}{\varphi}=u''$, so that $\nabla^2f=u''\cdot\mathrm{Id}_{TM}$, as well as $\Ric=-\frac{u^{(3)}}{(n_2-1)u'}\cdot\mathrm{Id}_{TM}$.
Therefore, $\nabla^2f=\frac{f}{n_2-1}\cdot\Ric$ if and only if
$u''=-\frac{uu^{(3)}}{(n_2-1)u'}$ on $I$.
But because $u''=\frac{\mu_1(f_1)}{2u}-\frac{n_2-2}{2u}(u')^2$, we have 
\begin{eqnarray*}
-\frac{uu^{(3)}}{(n_2-1)u'}&=&-\frac{u}{(n_2-1)u'}\cdot\left(\frac{\mu_1(f_1)}{
2u } -\frac{n_2-2}{2}\cdot\frac{(u')^2}{u}\right)'\\
&=&-\frac{u}{(n_2-1)u'}\cdot\left(-\frac{\mu_1(f_1)u'}{2u^2}-\frac{n_2-2}{2}
\cdot\frac{2u'u''u-(u')^3}{u^2}\right)\\
&=&\frac{1}{n_2-1}\cdot\left(\frac{\mu_1(f_1)}{2u}+\frac{n_2-2}{2}
\cdot\frac{2u''u-(u')^2}{u}\right)\\
&=&\frac{1}{n_2-1}\cdot\left(\frac{\mu_1(f_1)}{2u}+\frac{n_2-2}{2}
\cdot\frac{\mu_1(f_1)-(n_2-2)(u')^2-(u')^2}{u}\right)\\
&=&\frac{1}{n_2-1}\cdot\left(\frac{(n_2-1)\mu_1(f_1)}{2u}-\frac{
(n_2-2)(n_2-1)(u')^2}{2u}\right)\\
&=&u'',
\end{eqnarray*}
so that (\ref{eq:nabla2f1casef1=phin2>1}) is satisfied on $(M_1^2,g_1)$.\\
In the subcase where $n_2=2$, equation (\ref{eq:nabla2f2casef1=phi}) is 
equivalent to $(\nabla^{M_2})^2f_2=f_2\phi_2\cdot\mathrm{Id}_{TM_2}$, where 
$\phi_2:=\mu_1-\frac{S_2}{2}$.
Now (\ref{eq:Ricnablaf2n2>1}) 
yields $\frac{S_2}{2}\nabla^{M_2}f_2=\frac{S_2-\mu_1}{2}
\nabla^{M_2}f_2+\frac 
{f_2}{4}\nabla^{M_2}S_2$, 
that is, $f_2\nabla^{M_2}S_2=2\mu_1\nabla^{M_2}f_2$, which is equivalent to the 
existence of a real constant $C$ such that 
\[S_2=2\mu_1\ln(|f_2|)+C\]
on each connected component of the dense open subset $M_2\setminus 
f_2^{-1}(\{0\})$.
Denoting $\mu_2:=\mu_2(f_2)$, it can be deduced that
\begin{eqnarray*}
|\nabla^{M_2}f_2|^2&=&\frac{\mu_2}{2}-\frac{f_2^2(S_2-3\mu_1)}{2}\\
&=&\frac{\mu_2}{2}-\frac{f_2^2(2\mu_1\ln(|f_2|)+C-3\mu_1)}{2}\\
&=&\frac{\mu_2}{2}+\left(\frac{3\mu_1-C}{2}-\mu_1\ln(|f_2|)\right)f_2^2.
\end{eqnarray*}
This gives rise to a first-order ODE for $u(t):=f_2\circ F_t^\nu$, where 
$(F_t^\nu)_t$ is the local flow of 
$\nu:=\frac{\nabla^{M_2}f_2}{|\nabla^{M_2}f_2|_2}$ on some open subset of the
regular set of $f_2$.
Namely, \cite[Sec. 2]{Tashiro65} again implies that, on any open subset 
where $f_2$ has no critical point and vanishes nowhere, $(M_2^2,g_2)$ is 
locally isometric to 
$(\R^2,dt^2\oplus\rho(t)^2ds^2)$, where $\rho:=\frac{u'}{u'(0)}$.
Moreover, along any integral curve $\gamma$ of $\nu$, which is a geodesic of 
$(M_2,g_2)$ because of $\nabla^{M_2} f_2$ being a pointwise eigenvector of 
$(\nabla^{M_2})^2 f_2$, the function $u$ must satisfy the following 
first-order ODE: for any $t$ in some nonempty open interval,
\begin{equation}\label{eq:ODE1f2n2=2}
u'=\left(\frac{\mu_2}{2}+(\frac{3\mu_1-C}{2}
-\mu_1\ln(|u|))u^2\right)^{\frac{1}{2}}.
\end{equation}
Except in possibly very particular cases -- e.g. when $\mu_1=\mu_2=C=0$, in which $u$ is constant -- the maximal existence 
time for such a solution $u$ to (\ref{eq:ODE1f2n2=2}) is strictly contained in 
$\R$.
Note also that, if $u$ solves (\ref{eq:ODE1f2n2=2}), then 
\begin{eqnarray*}
u''&=&\frac{1}{2}\left(\frac{\mu_2}{2}+(\frac{3\mu_1-C}{2}
-\mu_1\ln(|u|))u^2\right)^{-\frac{1}{2}}\cdot\left((3\mu_1-C
-2\mu_1\ln(|u|))uu'-\mu_1u'u\right)\\
&=&\left(\frac{\mu_2}{2}+(\frac{3\mu_1-C}{2}
-\mu_1\ln(|u|))u^2\right)^{-\frac{1}{2}}\cdot\left(\mu_1-\frac{C}{2}
-\mu_1\ln(|u|)\right)uu'\\
&=&(u')^{-1}\cdot\left(\mu_1-\frac{C}{2}
-\mu_1\ln(|u|)\right)uu'\\
&=&\left(\mu_1-\frac{C}{2}
-\mu_1\ln(|u|)\right)u,
\end{eqnarray*}
where $\mu_1-\frac{C}{2}
-\mu_1\ln(|u|)=\mu_1-\frac{S_2\circ\gamma}{2}$ by the above identity for 
$S_2$.\\
This implies that, given any nowhere vanishing solution $u$ to 
(\ref{eq:ODE1f2n2=2}) on some open interval $I$ about $0$, the function 
$f(t,s):=u(t)$ solves 
\[\nabla^2f=u''\cdot\mathrm{Id}_{TM}=\left(\mu_1-\frac{S}{2}
\right)\cdot\mathrm{Id}_{TM}\]
on $(M_2^2,g_2):=(I\times\Sigma,dt^2\oplus(\frac{u'(t)}{u'(0)})^2ds^2)$, where 
$\Sigma=\R$ or $\mathbb{S}^1$.\\
Still in the case where $n_2=2$, equation (\ref{eq:nabla2f1casef1=phin2>1}) has 
not been considered yet in the literature as far as we know.
In the special subcase where $\mu_1=0$, which is equivalent to $S_1=0$, equation 
(\ref{eq:nabla2f1casef1=phin2>1}) can be rewritten under the form 
$(\nabla^{M_1})^2f_1=f_1\cdot\Ric_{M_1}-(\Delta^{M_1}f_1)\cdot\mathrm{Id}$, 
which is the general form of an element of $\ker(L_{g_1}^*)$ in 
\cite{Corvino00} 
when the underlying manifold is scalar-flat.
In case $\ker(L_{g_1}^*)\neq\{0\}$, the
metric $g_1$ is called \emph{static}. \color{black}
Although it is unclear whether a nonconstant positive solution $f_1$ to that 
equation can exist on a complete $M_1$, there is a noncomplete example: take the 
outer Schwarzschild manifold $(\R^3\setminus 
\overline{B}_{\frac{m}{2}},(1+\frac{m}{2r})^4\langle\cdot\,,\cdot\rangle)$ for 
some 
constant $m>0$, where $r=r(x)=|x|$ in $\R^3$ and 
$f_1(x)=\frac{1-\frac{m}{2r}}{1+\frac{m}{2r}}$, see \cite[p.145]{Corvino00}.
In case $M_1$ is either closed, complete with nonnegative Ricci curvature or 
with so-called moderate volume growth, the function $f_1$ must be constant.
The latter two are due to S.T.~Yau  \cite[Cor. 1 p. 217]{Yauharmonic75} and to L.~Karp 
\cite[Theorem B]{Karp82} (see also \cite[Sec. 3]{Karp82bis}) respectively, using only the harmonicity of $f_1$.
As a consequence, if $n_1=2$ (and $n_2=2$), then there is no nonconstant 
solution $f_1$ (for $S_1=0$ implies $\Ric_{M_1}=0$).\\
{\bf Case $n_2>2$ and $n_1>2$:} Then (\ref{eq:nabla2f1casef1=phin2>1}) defines 
a so-called $(0,n_1+n_2-1)$-Einstein metric on $(M_1,g_1)$ according to 
\cite{HePetersenWylie11102455,HePetersenWylie11102456} as we noticed in statement \ref{statement:f1}.
As for (\ref{eq:nabla2f2casef1=phi}), it has not been considered either in the 
literature when $\mu_1\neq0$ -- for $\mu_1=0$, it is already 
(\ref{eq:nabladf=-fric}) on $M_2$.
When $\mu_1\neq0$, we may take for $(M_2^{n_2},g_2,f_2)$ the standard solution 
to the Obata resp. Tashiro equation on the $n_2$-dimensional simply-connected 
spaceform of sectional curvature $\frac{\mu_1(f_1)}{n_2-2}$, which are the only 
Einstein solutions to (\ref{eq:nabla2f2casef1=phi}) when $n_2>2$.
This shows statement \ref{statement:exampleswp}.\\

In the particular case where $(M^n,g)$ is closed and $f\in
C^\infty(M,\R_+^\times)$ is such that $\mu_1(f):=k|\nabla f|^2-f\Delta f$ is
constant for some $k\in\R$, we can mimic the proof of Lemma \ref{l:eqnabladf=-fric}.\ref{claim:Mclosedfconstantsign}. 
First, we have $\mu_1(f)=0$: it suffices to evaluate $\mu_1(f)$
at two points, one where $\buil{\min}{M}(f)$ is attained and one where
$\buil{\max}{M}(f)$ is attained to obtain that $\mu_1(f)$ must be both
nonpositive and nonnegative because of $f>0$ and the opposite signs of the
Laplace operator of $f$ at a minimum and maximum respectively.
Independently, we can integrate $\mu_1(f)$ over $M$ and obtain
\[\mu_1(f)\cdot\mathrm{Vol}(M^n,g)=(k-1)\cdot\int_M|\nabla f|^2\, d\mu_g.\]
Therefore, if $k\neq1$, then $f$ must be constant.
If $k=1$, the vanishing of $\mu_1(f)$ is equivalent to $\Delta f=\frac{|\nabla f|^2}{f}\geq0$ on the closed
manifold $M$, which with $\int_{M}\Delta f\,d\mu_g=0$ shows
that, again, $\nabla f=0$ must hold on $M$, therefore $f$ must also
be constant on $M$.
This proves statement \ref{statement:Mclosed} and concludes the proof of Proposition \ref{p:hessianricciwarpedproducts}.
\findemo

In case the factor $(M_1,g_1)$ of the warped product is complete, we show that actually the map $f$ must be constant along $M_1$.

\begin{prop}\label{p:mumu1mu2}

Let $f=\pi_1^*f_1\cdot\pi_2^*f_2$ satisfy {\rm(\ref{eq:nabladf=-fric})} on $(M^n,g)=(M_1\times M_2, g_1\oplus f_1^2 g_2)$  for some smooth positive function $f_1$ on $M_1$ and smooth function $f_2$ on $M_2$.
Assume $(M_1,g_1)$ to be complete and connected.\\
Then $f_1$ must be constant on $M_1$, the manifold $(M_1,g_1)$ must be Ricci-flat and $f_2$ must satisfy {\rm(\ref{eq:nabladf=-fric})} on $(M_2,g_2)$.
Therefore, the map $W(M_2,g_2)\longrightarrow W(M,g)$ extending any solution {\rm(\ref{eq:nabladf=-fric})} to $M$ is an isomorphism.
\end{prop}

{\it Proof:} In case $f_1>0$ on $M_1$ and for $f=\pi_1^*f_1\cdot\pi_2^*f_2$ on $M_1\times_{f_1^2} M_2$, the quantities $\mu(f)$, $\mu_1(f_1)$ and $\mu_2(f_2)$ defined above are related as follows:
\begin{eqnarray*}
\mu(f)&=&f\Delta f+2|\nabla f|^2\\
&=&f_1f_2((\Delta f_1)f_2+f_1\Delta f_2)+2|f_2(\nabla f_1)+f_1\nabla f_2|^2\\
&=&f_1f_2((\Delta^{M_1} f_1)f_2+\frac{f_1}{f_1^2}\Delta^{M_2} f_2)+2|f_2(\nabla^{M_1} f_1)+\frac{f_1}{f_1^2}\nabla^{M_2} f_2|^2\\
&=&f_1(\Delta^{M_1}f_1)f_2^2+f_2(\Delta^{M_2}f_2)+2f_2^2|\nabla^{M_1}f_1|_1^2+2|\nabla^{M_2}f_2|_2^2\\
&=&\left(f_1(\Delta^{M_1}f_1)+2|\nabla^{M_1}f_1|_1^2\right)\cdot f_2^2+f_2\Delta^{M_2}f_2+2|\nabla^{M_2}f_2|_2^2\\
&=&\left(f_1(\Delta^{M_1}f_1)+2|\nabla^{M_1}f_1|_1^2+\mu_1(f_1)\right)\cdot f_2^2\\
&&+f_2\Delta^{M_2}f_2+2|\nabla^{M_2}f_2|_2^2-\mu_1(f_1)f_2^2\\
&=&n_2|\nabla^{M_1}f_1|_1^2f_2^2+\mu_2(f_2).
\end{eqnarray*}
This implies that, if $f\neq0$ solves (\ref{eq:nabladf=-fric}) and $\varphi=f_1>0$, then $|\nabla^{M_1} f_1|_1$ is constant on $M_1$.
Note that this holds whether $(M_1,g_1)$ is complete or not, i.e. whenever $M_1$ is connected.
From now on assume $(M_1,g_1)$ to be complete.
By contradiction, if $|\nabla^{M_1} f_1|_1$ were a positive constant, then $f_1$ would have no critical point on $M_1$ and therefore the flow of the normalised gradient vector field $\nu_1:=\frac{\nabla^{M_1} f_1}{|\nabla^{M_1} f_1|_1}$ would define a diffeomorphism from $M_1$ to the product $\R\times\Sigma_1$ for some smooth level hypersurface $\Sigma_1$ of $f_1$; and $f_1$ would be a nonconstant affine linear function of $t\in\R$.
But this would contradict $f_1>0$ on $M_1$.
Therefore, $\nabla^{M_1} f_1=0$ must hold on $M_1$ i.e., $f_1$ must be constant on $M_1$.
In turn, this implies that $\mu_1(f_1)=0$, $\Ric_{M_1}=0$ when $n_2\geq2$ (anyway $\Ric_{M_1}=0$ when $n_2=1$ as we saw above) and that $f_2\in W(M_2,g_2)$.
Therefore, the function $f$ is the trivial extension on $M$ of $f_2\in W(M_2,g_2)$.
\findemo
\color{black}

\section{Case where $\dim(W(M^n,g))\geq2$}\label{s:dimWgeq2}

In this section, we look at the particular case
where (\ref{eq:nabladf=-fric}) has a $k\geq2$-dimensional space of solutions 
and then look at \emph{homogeneous} examples.
\begin{prop}\label{p:dimWgeq2}
Let $(M^n,g)$ be any connected complete Riemannian manifold.
Assume that {\rm (\ref{eq:nabladf=-fric})} has a $k\geq2$-dimensional space of 
solutions.
Then we have one of the following:

\begin{enumerate}
\item Case $k=2$: the manifold $(M^n,g)$ must be isometric to the Riemannian product $(M_1^{n-1}\times\R,g_1\oplus dt^2)$ for some complete Ricci-flat manifold admitting no line  $(M_1^{n-1},g_1)$.
Moreover, the solutions of {\rm(\ref{eq:nabladf=-fric})} on $(M^n,g)$ are the affine linear functions of $t\in\R$ extended constantly along $M_1$.
\item Case $k>2$: the manifold $(M^n,g)$ must be isometric to the Riemannian product $(M_1^{n-k+1}\times M_2^{k-1},g_1\oplus g_2)$ for some complete Ricci-flat manifold admitting no line  $(M_1^{n-k+1},g_1)$ and where $(M_2^{k-1},g_2)$ is either $\mathbb{S}^2,\R^2$ or $\mathbb{H}^2$ with standard metric of curvature $1,0,-1$ (up to rescaling $g$) for $k=3$ or is $\R^{k-1}$ with standard flat metric for $k>3$.
Moreover, the solutions of {\rm(\ref{eq:nabladf=-fric})} on $(M^n,g)$ are the solutions of the Obata resp. Tashiro equation on $(M_2,g_2)$ extended constantly along $M_1$.
\end{enumerate}
%
\color{black}
\end{prop}

{\it Proof:} We first assume $M$ to be simply-connected.
By \cite[Theorem A]{HePetersenWylie11102455}, which can be applied since 
(\ref{eq:nabladf=-fric}) is the particular case of the equation 
$\nabla^2f=f\cdot q$ for some quadratic form $q$ on $TM$, we already know that, 
if $k\geq2$, then $(M^n,g)$ must be isometric to the warped product $(M_1\times 
M_2,g_1\oplus f_1^2 g_2)$ for some smooth positive function $f_1$ on 
$M_1$, where $(M_1^{n-k+1},g_1)$ and $(M_2^{k-1},g_2)$ are complete \cite[Lemma 
7.2]{BishopONeill69} simply-connected 
Riemannian manifolds and $f_1$ is a smooth positive function on $M_1$.
Moreover, $(M_2,g_2)$ must be a spaceform and any solution $f$ of 
(\ref{eq:nabladf=-fric}) is of the form $f=\pi_1^*f_1\cdot\pi_2^*f_2$, where 
$f_2$ satisfies the Obata resp. Tashiro equation on $(M_2,g_2)$ 
\cite[Theorem B]{HePetersenWylie11102455}.
Taking the above considerations on solutions of (\ref{eq:nabladf=-fric}) on 
warped products into account in case $f_1$ is the warping function, Proposition \ref{p:mumu1mu2} can be applied and implies that $f_1$ is constant, that $(M_1,g_1)$ is Ricci-flat and that $f_2\in W(M_2,g_2)$.
We look at different cases according to $k$:
\begin{enumerate}
\item Case $k=2$: then we could conclude above that $f_2$ is an affine linear function of $t\in\R$.
Since no nonconstant affine function can be periodic, any group action leaving invariant some nonconstant $f_2\in W(M_2,g_2)$ must be trivial.
Moreover, if $(M_1,g_1)$ could be split off a line, then it would be isometric to $\Sigma_1\times\R$ for some smooth hypersurface $\Sigma_1$ of $M_1$; but then $M_1\times\R\cong\Sigma_1\times\R^2$ would carry a $k\geq3$-dimensional space of solutions to (\ref{eq:nabladf=-fric}), which would contradict $k=2$.
Therefore, $(M_1,g_1)$ cannot contain any line.
\item Case $k>2$: then we could conclude above that $f_2\in W(M_2,g_2)$.
If $k=3$, then, up to rescaling $g$, the manifold $(M_2,g_2)$ must be isometric to either $\mathbb{S}^2,\R^2$ or $\mathbb{H}^2$ with standard metric of constant curvature $1,0,-1$ respectively; and $W(M_2,g_2)$
must consist of the solutions of the Obata resp. Tashiro equation on $(M_2,g_2)$ as we saw in Lemma \ref{l:eqnabladf=-fric}.\ref{claim:MEinstein}.
Again, in case $M_2=\mathbb{S}^2$ or $\mathbb{H}^2$, no group action on $M_2$ can leave any nonzero solution to (\ref{eq:nabladf=-fric}) invariant on $M_2$.
If $M_2=\R^2$, then no nontrivial group action preserves the $3$-dimensional space of affine linear functions on $\R^2$.\\
If $k>3$, then, as a consequence of Lemma \ref{l:eqnabladf=-fric}.\ref{claim:MEinstein}, the manifold $(M_2,g_2)$ must be isometric to flat $\R^{k-1}$ and again no nontrivial group action preserves the $k$-dimensional space of affine linear functions on $\R^{k-1}$.\\
In both subcases, $(M_1,g_1)$ cannot contain any line, otherwise $\dim(W(M^n,g))\geq k+1$.
\end{enumerate}
In all cases, the only possible nontrivial group actions on $M_1\times M_2$ is trivial along the $M_2$ factor.
Thus, if $M$ is not simply-connected, then $M$ must be isometric to $M_1^{n-k+1}\times M_2^{k-1}$, where $M_2$ is a simply connected model space as above and $M_1$ is a complete Ricci-flat manifold having no line since its universal cover cannot contain any.
Furthermore, every $f\in W(M,g)$ must be the trivial extension on $M_1\times M_2$ of a solution $f_2\in W(M_2,g_2)$.
This concludes the proof of Proposition \ref{p:dimWgeq2}.

\findemo

Note that, as a consequence of Proposition \ref{p:dimWgeq2}, if a complete $(M^n,g)$ carries an $(n+1)$-dimensional space of solutions to (\ref{eq:nabladf=-fric}) with $n\neq2$, then $(M^n,g)$ must be isometric to $\R^{n+1}$ with standard flat metric.
\color{black}


\section{Homogeneous case}\label{s:homogcase}

Next, we look at homogeneous manifolds carrying nontrivial solutions of 
(\ref{eq:nabladf=-fric}).

\begin{prop}\label{p:Mhomogeneous}
Let $(M^n,g)$ be any connected
homogeneous Riemannian
manifold.
Assume the existence of a non-iden\-ti\-cal\-ly vanishing smooth function $f$ 
on $M$ 
satisfying {\rm (\ref{eq:nabladf=-fric})}.\\
Then one of the following holds:

\begin{enumerate}
\item\label{statement:scal0homogeneous} If the scalar curvature $S$ of $(M^n,g)$ vanishes and $f$ is nonconstant, then $(M^n,g)$ must be isometric to a flat manifold $\rquot{\R^n}{\Gamma}$ for some discrete fixed-point free subgroup $\Gamma$ of $\mathrm{O}(n)\ltimes\R^n$.
\item\label{statement:k=2homogeneous} If $k:=\dim(W(M^n,g))=2$, then $(M^n,g)$ must be isometric to the Riemannian product $\rquot{\R^{n-1}}{\Gamma}\times\R$ for some discrete fixed-point free and co-compact subgroup $\Gamma$ of $\mathrm{O}(n-1)\ltimes\R^{n-1}$.
In that case, the map $W(\R,dt^2)\longrightarrow W(M^n,g)$ extending any affine linear function trivially on the  first factor is an isomorphism.
\item\label{statement:k=3homogeneous} If $k=3$, then up to re\-sca\-ling $g$, the manifold $(M^n,g)$ must be
isometric to the Riemannian product $\rquot{\R^{n-2}}{\Gamma}\times S^2(\varepsilon)$, where
$S^2(\varepsilon)$ is the simply-connected complete surface of constant curvature
$\varepsilon\in\{0,\pm1\}$ and $\rquot{\R^{n-2}}{\Gamma}$ is a compact flat manifold.
In that case, the map $W(S^2(\varepsilon),g_{S^2(\varepsilon)})\longrightarrow W(M^n,g)$ extending any function trivially on the $\Sigma$-factor is an isomorphism.
\item\label{statement:k>=4homogeneous} If $k\geq4$, then $(M^n,g)$ must be isometric to the Riemannian product $\rquot{\R^{n-k+1}}{\Gamma}\times\R^{k-1}$, where $\rquot{\R^{n-k+1}}{\Gamma}$ is a compact flat manifold and $\R^{k-1}$ carries its standard Euclidean metric.
\item\label{statement:k=1homogeneous} If $k=1$,
then unless $W(M^n,g)$ consists of constant functions,
$\mu(f)=0$ must hold for every $f\in W$.
 Moreover, the manifold $(M^n,g)$ must be a one-dimensional
extension of some homogeneous Riemannian manifold satisfying the particular
conditions {\rm (\ref{eq:SARichomog})} below.
\end{enumerate}

\end{prop}

{\it Proof:} If $(M^n,g)$ has vanishing scalar curvature and $f$ is nonconstant, then we already know from Lemma \ref{l:eqnabladf=-fric} that $(M^n,g)$ must be Ricci-flat.
But because any homogeneous Ricci-flat Riemannian manifold must be flat  \cite{AlekseevskiiKimelfeld75}, actually $(M^n,g)$ must be isometric to a flat manifold $\rquot{\R^n}{\Gamma}$ for some discrete and necessarily fixed-point free subgroup $\Gamma$ of $\mathrm{O}(n)\ltimes\R^n$.
This shows statement \ref{statement:scal0homogeneous}.\\
If $\dim(W(M^n,g))=k\geq2$, then Proposition \ref{p:dimWgeq2} implies that $(M^n,g)$ must be isometric to the Riemannian product $M_1^{n-k+1}\times M_2^{k-1}$, where $M_1^{n-k+1}$ is a Ricci-flat manifold containing no line and $M_2^{k-1}$ is flat Euclidean space except when $k=3$, in which case it is also allowed to be $\mathbb{S}^2$ or $\mathbb{H}^2$ with standard spherical resp. hyperbolic metric.
Moreover, any solution to ({\ref{eq:nabladf=-fric}) must be the trivial extension to $M$ of a standard solution on $M_2$.
Now recall the following result, which is a combination of Lemma 5.6 and the first part of the proof of
Theorem 5.7 in \cite{HePetersenWylie11102456}; the latter can be applied
because of $W(M^n,g)$ being invariant under isometry: in our notation, the
isometries of $(M_1\times M_2,g_1\oplus g_2)$ are the maps of the form
$h=(h_1,h_2)$, where $h_1$ and $h_2$ are isometries of $(M_1,g_1)$ and $(M_2,g_2)$ respectively.
This already implies that, writing $M=\rquot{G}{K}$, the group $G$ when
 can be embedded into the direct product of two groups, the first one
acting isometrically and transitively on $M_1$ and the second one acting
transitively on $M_2$.
In particular, $(M_1,g_1)$ must itself be homogeneous.
In turn, this implies that, being Ricci-flat, $(M_1,g_1)$ must be flat, again by \cite{AlekseevskiiKimelfeld75}.
Therefore $(M_1,g_1)$ must be isometric to $\rquot{\R^{n-k+1}}{\Gamma}$ for some discrete fixed-point free subgroup $\Gamma$ of $\mathrm{O}(n-k+1)\ltimes\R^{n-k+1}$.
Since only \emph{compact} flat manifolds have no line, the subgroup $\Gamma$ must be co-compact i.e., $M_1$ must be compact.
This shows statements \ref{statement:k=2homogeneous}, \ref{statement:k=3homogeneous} and \ref{statement:k>=4homogeneous}.
\color{black}

Let us now assume the space $W(M^n,g)$ of functions satisfying 
(\ref{eq:nabladf=-fric}) to be one-dimensional on $M=\rquot{G}{K}$.
Then as in 
\cite[Sec. 5]{HePetersenWylie11102456} we consider the action of $G$ on $W(M^n,g)$.
Because the Ricci-tensor of $M$ is isometry- and thus $G$-invariant, so is 
equation (\ref{eq:nabladf=-fric}), i.e. for every $f$ satisfying 
(\ref{eq:nabladf=-fric}) and every $h\in G$, the function $f\circ L_{h^{-1}}$ 
also satisfies (\ref{eq:nabladf=-fric}).
But because of $\dim(W(M^n,g))=1$, there exists for a fixed nonzero $f\in 
W(M^n,g)$ and 
every $h\in G$ a nonzero constant $C_h$ such that $f\circ L_{h^{-1}}=C_h\cdot 
f$.
The map $G\to\R^\times$, $h\mapsto C_h$ is a Lie-group homomorphism and 
actually takes its values in $\{\pm1\}$ if $\mu(f)\neq0$ since, by invariance 
of $\mu(f)$ under isometry,
\[\mu(f)=\mu(f\circ L_{h^{-1}})=\mu(C_h\cdot f)=C_h^2\cdot\mu(f)\]
for every $h\in G$.
Therefore, if $\mu(f)\neq0$, then $C_h\in\{\pm1\}$ for every $h\in G$.
Now if $M$ is connected as in the assumptions, then so can be assumed $G$ 
(otherwise replace $G$ by the connected component of the neutral element), in 
which case necessarily $C_h=1$ holds for every $h\in G$ and therefore every 
$f\in W(M^n,g)$ is constant.

Therefore $\mu(f)=0$ holds.
As a consequence, $S=-2$ and $f$ has no critical point on $M$, see Lemma 
\ref{l:eqnabladf=-fric}.\\
Next we show that $(M^n,g)$ must be the one-dimensional extension of some 
homogeneous Riemannian manifold $N^{n-1}$ with Ricci-tensor having particular 
properties.
Consider the subgroup $H$ of $G$ defined by
\[H:=\left\{h\in G\,|\,C_h=1\right\},\]
that is, $H$ is the subgroup of all elements of $G$ leaving a (thus any) 
function $f\in W(M^n,g)$ invariant.
Since $C\colon G\to\R_+^\times$ is a nontrivial and therefore surjective 
Lie-group-homomorphism, $H=\ker(C)$ is a closed normal subgroup of $G$ and of 
codimension $1$.
Moreover, fixing $f\in W(M^n,g)\setminus\{0\}$, we know from Lemma 
\ref{l:eqnabladf=-fric} that $f(M)=\R_+^\times=(0,\infty)$ since $f$ can 
be expressed as an exponential function along any integral curve of its 
normalised gradient.
We let $N:=f^{-1}(\{1\})$, which is a smooth hypersurface of $M$.
By definition, $H$ leaves $N$ invariant.
Moreover, fixing some $x\in N$, any $h\in G$ with $L_h(x)=x$ must satisfy 
$C_h=1$ and therefore lie in $H$.
In other words, the isotropy group $H_x:=\left\{h\in 
H\,|\,L_h(x)=x\right\}$ of $x$ under the $H$-action must coincide with $K=G_x$.
Independently, for any $y\in N$, there is an $h\in G$ such that $L_h(x)=y$; 
again, because of $f(x)=f(y)\neq0$, necessarily $C_h=1$ must hold, i.e. $h\in 
H$.
This proves that the orbit $H\cdot x:=\left\{L_h(x)\,|\,h\in H\right\}$ of $x$ 
in $N$ must be all of $N$ and therefore $N=\rquot{H}{K}$ is a $H$-homogeneous 
Riemannian manifold.
As in the proof of \cite[Theorem 5.1]{HePetersenWyliehomogeneous2015}, we 
split the Lie algebra $\underline{G}=\underline{P}\oplus\underline{K}$ of $G$ 
in an $\mathrm{Ad}_G(K)$-invariant and orthogonal way and let 
$\xi\in\underline{P}\cong TM$ be the vector corresponding to $\nu\in T^\perp N$.
Note that, because of $C_{|_H}=1$, the gradient of $f$ and therefore also $\nu$ 
are preserved by the $H$-action, so that $\xi$ makes sense.
Actually, 
$\underline{P}=\R\xi\oplus\left((\R\xi)^\perp\cap\underline{P}\right)$ and 
$\underline{H}=\left((\R\xi)^\perp\cap\underline{P}\right)\oplus\underline{K}$, 
the splittings being orthogonal.
Furthermore, the Lie-bracket of $\xi$ in $\underline{G}$ preserves 
$\underline{H}$ because of $H$ being a normal subgroup of $G$.
This already proves that $G=H\ltimes\R$ and that $(M,g)$ is
the one-dimensional extension of the $H$-homogeneous space 
$(N^{n-1},g_{|_N})$.\\
In that case, following \cite{HePetersenWyliehomogeneous2015}, we fix some 
$\alpha\in\R^\times$ and let $D:=\frac{1}{\alpha}[\xi,\cdot]=\frac{1}{\alpha}\mathcal{L}_\xi$, which is hence a 
derivation of $\underline{H}$.
We denote by $\mathcal{S}$ and $\mathcal{A}$ the symmetric and skew-symmetric 
parts of $D$ respectively seen as endomorphisms of $TN$, see 
\cite[Eq. (2.1)]{HePetersenWyliehomogeneous2015}.
Let $\mathcal{T}:=-\nabla\xi$ denote the Weingarten map of $N$ in $M$.
Then by \cite[Prop. 2.7]{HePetersenWyliehomogeneous2015} we have 
$\mathcal{T}=\alpha\mathcal{S}$ and 
$\nabla_\xi\mathcal{T}=-\alpha^2[\mathcal{S},\mathcal{A}]$.
Furthermore, \cite[Lemma 2.9]{HePetersenWyliehomogeneous2015} implies that, for 
all $X,Y\in TN$,
\[\left\{\begin{array}{ll}\mathrm{ric}(\xi,\xi)&=-\alpha^2\mathrm{tr}
(\mathcal{S}^2)\\\mathrm{ric}(X,\xi)&=\alpha(\delta\mathcal{S})(X)\\\mathrm{ric}
(X,Y)&=\mathrm{ric}^N(X,Y)-(\alpha^2\mathrm{tr}(\mathcal{S}))g(\mathcal{
S}X,Y)-\alpha^2g([\mathcal{S},\mathcal{A}]X,Y)\end{array}\right.\]
Now writing $f(t)=e^t$, where $t$ lies in the $\R$-factor of $G=H\ltimes\R$, we
have $\nabla df=fdt^2-fg(T\cdot,\cdot)$ which, together with $\nabla_\xi\xi=0$, 
gives that identity (\ref{eq:nabladf=-fric}) is equivalent to 
\[
\left\{\begin{array}{ll}
\alpha^2\mathrm{tr}
(\mathcal{S}^2)(=\alpha^2|\mathcal{S}|^2)&=1\\
\alpha(\delta\mathcal{S})&=0\\
-\alpha 
g(\mathcal{S}X,Y)&=-\mathrm{ric}^N(X,Y)+\alpha^2\mathrm{tr}(\mathcal{S}
)g(\mathcal{S}X,Y)+\alpha^2g([\mathcal{S},\mathcal{A}]X,Y)
\end{array}
\right.
\]
for all $X,Y\in TN$.
In other words, (\ref{eq:nabladf=-fric}) is equivalent to 
\begin{equation}\label{eq:SARichomog}
\left\{\begin{array}{ll}\alpha&=\frac{\epsilon}{|\mathcal{S}|}\\
\delta\mathcal{S}&=0\\
\mathrm{Ric}_N&=\frac{1}{|\mathcal{S}|^2}\left(\left(\mathrm{tr}(\mathcal{S})
+\epsilon|\mathcal{S}|\right)\mathcal{S}+[\mathcal{S},\mathcal{A}]\right)
\end{array}
\right.
\end{equation}
for some $\epsilon\in\{\pm1\}$.
This shows statement \ref{statement:k=1homogeneous} and completes the proof of Proposition \ref{p:Mhomogeneous}.
\findemo

The case where $\dim(W(M^n,g))=1$ could 
lead to new examples, see \cite{HePetersenWyliehomogeneous2015} and 
\cite{HePetersenWyliecsc}.\\

\color{black}

\section{K\"ahler case}\label{s:Kaehlercase}

As in \cite{CaseShuWei2011}, we next consider the case where $(M^n,g)$ is 
assumed to be K\"ahler:

\begin{prop}
Assume $(M^{2n},g,J)$ to be a complete K\"ahler manifold and let $f$ be any 
nonconstant smooth real-valued function satisfying {\rm 
(\ref{eq:nabladf=-fric})} on $M$.
Then, up to re\-sca\-ling $g$, the K\"ahler manifold $(M^{2n},g,J)$ 
is holomorphically isometric to $S^2(\varepsilon)\times\Sigma^{2n-2}$ for some 
Ricci-flat K\"ahler manifold $\Sigma$, where $S^2(\varepsilon)=\mathbb{S}^2$ if 
$\varepsilon=1$, $\mathbb{H}^2$ if $\varepsilon=-1$ and either $\R^2$ or 
$\R\times\mathbb{S}^1$ if $\varepsilon=0$; moreover, the K\"ahler structure is the 
product K\"ahler structure and $f$ is the trivial extension to $M$ of a solution 
to {\rm 
(\ref{eq:nabladf=-fric})} on $S^2(\varepsilon)$.
\end{prop}

{\it Proof:} The first steps follow those in the proof of 
\cite[Theorem 1.3]{CaseShuWei2011}.
Since the Ricci-tensor of $(M,g,J)$ is $J$-invariant, so is the Hessian of $f$ 
by (\ref{eq:nabladf=-fric}), i.e. $\nabla^2 f\circ J=J\circ\nabla^2 f$.
As a first consequence, the vector field $J\nabla f$ is a (real) holomorphic 
vector field on $(M,g,J)$ and therefore its zeros -- which are precisely the 
critical points of $f$ -- form a totally geodesic K\"ahler submanifold of $M$ 
of dimension $2k<2n$; in particular the regular set of $f$ is dense in $M$.
As a se\-cond consequence, the $2$-form $g(\nabla^2f\circ J\cdot\,,\cdot)$ may 
be rewritten $\frac{1}{2}\mathcal{L}_{\nabla f}\Omega$, where $\Omega:=g(J\cdot\,,\cdot)$ 
is the K\"ahler form of $(M,g,J)$.
Therefore, 
\[d\left(g(\nabla^2f\circ J\cdot\,,\cdot)\right)=\frac{1}{2}d\left(\mathcal{L}_{\nabla 
f}\Omega\right)=\frac{1}{2}d\left(\nabla f\lrcorner d\Omega+d(\nabla 
f\lrcorner\Omega)\right)=0,\]
i.e. $g(\nabla^2 f\circ J\cdot\,,\cdot)$ is a closed $2$-form on $M$.
But because the Ricci-form $g(\Ric\circ J\cdot\,,\cdot)$ is also closed on $M$, 
so is the $2$-form $\frac{1}{f}g(\nabla^2 f\circ J\cdot\,,\cdot)$ on 
$\{f\neq0\}$, again by (\ref{eq:nabladf=-fric}).
This implies $df\wedge \left(g(\nabla^2 f\circ J\cdot\,,\cdot)\right)=0$ on 
$\{f\neq0\}$ and therefore on $M$ by density (recall that $f^{-1}(\{0\})$, if 
nonempty, is a totally geodesic hypersurface of $(M,g)$).
In turn this implies the existence at each regular point of $f$ of a linear 
form $\lambda$ on $(\nabla f)^\perp$ such that, for every $X\perp\nabla f$,
\begin{equation}\label{eq:nabla2JXf}
\nabla_{JX}^2f=\lambda(X)\nabla f.
\end{equation}
For $X=J\nabla f$, we obtain via (\ref{eq:Ricnablaf}) that $\nabla S$ is 
pointwise tangent to $\nabla f$, i.e. there exists a function $\theta$ on $M$ 
such that $\nabla S=\theta\nabla f$ on $M$ (this holds true on the regular set 
of $M$ and hence on $M$ by density, taking into account that at every critical 
point both $\nabla f$ and $\nabla S$ vanish).
For $X\in\{\nabla f,J\nabla f\}^\perp$, by $J$-invariance of $\nabla^2f$ the 
r.h.s. of (\ref{eq:nabla2JXf}) must vanish whenever the basepoint is a regular 
point of $f$.
In turn this implies $\Ric(X)=0$ for all $X\in\{\nabla f,J\nabla f\}^\perp$ and 
at every regular point of $f$.
Now because of $\Ric(\nabla f)=\left(\frac{S}{2}+\frac{f\theta}{4}\right)\nabla 
f$, the $J$-invariance of $\Ric$ and $\Ric_{|_{\{\nabla f,J\nabla 
f\}^\perp}}=0$, we obtain
\[S=S+\frac{f\theta}{2},\]
so that $\theta=0$, first on the regular set of $f$ and then on $M$ by density, 
i.e. $S$ is constant on $M$.
This implies that both distributions $\mathrm{Span}(\nabla f,J\nabla f)$ and 
$\{\nabla f,J\nabla f\}^\perp$ are integrable and totally geodesic, the former 
one being the tangent bundle of a surface of curvature $\frac{S}{2}$ -- which 
may be assumed to be $\pm 1$ up to re\-sca\-ling $g$ in case $S\neq0$ -- and the 
latter the tangent bundle of a necessarily Ricci-flat K\"ahler manifold 
$\Sigma$.
The rest of the proof is analogous to that of Proposition \ref{l:eqnabladf=-fricscalconstant}.\ref{statement:scalcstn=3}.
\findemo

\end{document}